\documentclass{amsart}

\usepackage{amsmath,amssymb,amsthm,texdraw,citesort}

\theoremstyle{plain}
\newtheorem{theorem}{Theorem}[section]
\newtheorem{proposition}[theorem]{Proposition}
\newtheorem{lemma}[theorem]{Lemma}
\theoremstyle{definition}
\newtheorem{definition}[theorem]{Definition}
\theoremstyle{remark}
\newtheorem{example}[theorem]{Example}

\newcommand{\ldot}{%
\raisebox{-0.195em}{\rule{1pt}{0pt}$\cdot$\rule{1pt}{0pt}}}
\newcommand{\la}{\lambda}

\newcommand{\calT}{\mathcal{T}}
\newcommand{\rp}{\hat{P}^+}
\newcommand{\TL}{\calT^L}
\newcommand{\Tinf}{\calT(\infty)}
\newcommand{\fit}{\tilde{f}_i}

\newcommand{\eit}{\tilde{e}_i}
\newcommand{\veps}{\varepsilon}
\newcommand{\vphi}{\varphi}
\newcommand{\Binf}{\mathcal{B}(\infty)}
\newcommand{\bpl}{\bar{\pi}_{\la}}
\newcommand{\iso}{\psi}
\newcommand{\Z}{\mathbf{Z}}
\newcommand{\B}{\mathcal{B}}
\newcommand{\ol}{\overline}

\DeclareMathOperator{\wt}{wt}

\begin{document}

\title
[Young tableaux and crystal $\Binf$]
{Young tableaux and crystal $\Binf$\\
for finite simple Lie algebras}

\author{Jin Hong}
\address{National Security Research Institute\\
         161 Gajeong-dong, Yuseong-gu\\
         Daejeon, 305-350, Korea}
\email{jinhong@etri\ldot re\ldot kr}
\author[Hyeonmi Lee]{Hyeonmi Lee$^*$}
\address{Korea Institute for Advanced Study\\
         207-43 Cheongnyangni 2-dong, Dongdaemun-gu\\
         Seoul 130-722, Korea}
\email{hmlee@kias\ldot re\ldot kr}
\thanks{$^*$This work was supported in part by KOSEF Grant
        R01-2003-000-10012-0}
\date{July 22, 2005}

\begin{abstract}
We study the crystal base of the negative part of a quantum group. An
explicit realization of the crystal is given in terms of Young tableaux
for types $A_n$, $B_n$, $C_n$, $D_n$, and $G_2$. Connection between our
realization and a previous realization of Cliff is also given.
\end{abstract}

\maketitle

\section{Introduction}

Quantum group $U_q(\mathfrak{g})$ is a $q$-deformation of the universal
enveloping algebra over a Lie algebra $\mathfrak{g}$, and crystal bases
reveal the structure of $U_q(\mathfrak{g})$-modules in a very
simplified form. As these $U_q(\mathfrak{g})$-modules are known to be
$q$-deformations of modules over the original Lie algebras, knowledge
of these structures also affects the study of Lie algebras.

The crystal $\Binf$, which is the crystal base of the negative part
$U_q^-(\mathfrak{g})$ of a quantum group, has received attention since
the very birth of crystal base theory~\cite{K91,K90}. This is not only
because it is an essential part of the \emph{grand loop} argument
proving the existence of crystal bases, but because it gives insight
into the structure of quantum group itself.

Much effort has been made~\cite{MR1475048,saito,cliff,lee,hm} to give
explicit description of the crystals $\Binf$ over various Kac-Moody
algebras. A related well known result is that it is possible to
describe the highest weight crystal $\B(\la)$ over finite simple Lie
algebras, in terms of Young tableaux~\cite{KN,KM}. This has lead to the
belief that it should be possible to give a similar description for
$\Binf$ also. So far, only the $A_n$ type has been dealt with in this
direction.

In the current work, we restrict ourselves to finite simple Lie
algebras of types $A_n$, $B_n$, $C_n$, $D_n$, and $G_2$. These are the
cases for which Young tableau realization of crystal $\B(\la)$ is
known. For these cases, we explicitly describe $\Binf$ in terms of
Young tableaux. Our result for $A_n$ type will be equivalent to that
of~\cite{hm}, but it will be obtained through a completely different
approach. As most of these cases were dealt with in
Cliff's~\cite{cliff} realization given in terms of a completely
different object, we give an isomorphism between our realization and
that of Cliff.

The previous work~\cite{hm} gives a Young tableaux realization of
$\Binf$ for the $A_n$ case, and then uses this to create yet another
realization based on Nakajima monomials, an object which was originally
introduced to characterize highest weight
crystals~\cite{MR1988989,MR1988990}. We expect results of this paper
also to lead to similar realizations of $\Binf$ in terms of Nakajima
monomials for the remaining finite types $B_n$, $C_n$, $D_n$, and
$G_2$.

The paper is organized as follows. We start by introducing the notion
of large semi-standard tableaux. Then, in the following section, these
are collected together into a set, an equivalence relation is given to
the set, and a crystal structure is given to the resulting set of
equivalence classes. In Section~\ref{sec:4}, this new crystal is shown
to be isomorphic to $\Binf$. Our main result is given in
Section~\ref{sec:5}, where a set of representatives for our new crystal
is explicitly presented. This gives a new realization of $\Binf$. The
last section connects our result with another realization of $\Binf$
given by Cliff~\cite{cliff}.

\section{Large semi-standard tableaux}

Throughout this paper, we shall be dealing with finite Lie algebras of
types $A_n$, $B_n$, $C_n$, $D_{n+1}$, and $G_2$. Unless explicitly
stated otherwise, all our discussions will hold true for each of these
types. Notice that the subscript for $D$-type is different from the
others. This is to simplify our later writings, and does not imply any
restriction on the range of $D$-types we are considering. For the $G_2$
case, $n=2$ should be assumed.

We shall assume knowledge of the basic theory of crystal bases, and
related standard notation, for example, as given in the
books~\cite{MR1881971,MR1359532}, will be used.

The crystal base of $U_q^-(\mathfrak{g})$, first introduced
in~\cite{K91}, will be denoted by $\Binf$. For each of the finite
classical types, we shall use the definitions of \emph{semi-standard}
tableaux as given by Kashiwara and Nakashima~\cite{KN}. For the $G_2$
type, we shall take the Young tableau realization of highest weight
crystal $\B(\la)$ given in~\cite{KM} as the definition of semi-standard
tableaux. Since the first of these two works is a rather well known
result, and since the second is very similar in spirit to the first, we
refer readers to the original papers and shall not repeat the
complicated definitions here. The alphabet to be used inside the boxes
constituting the Young tableaux for each type will be denoted commonly
by $J$, and it will be equipped with an ordering $\prec$, as given
in~\cite{KN,KM}. For example, in the $C_n$ case, it would be
\begin{equation*}
J = \{ 1 \prec 2 \prec \cdots \prec n \prec \bar{n} \prec \cdots \prec
\bar{2} \prec \bar{1} \}.
\end{equation*}
Also, based on results of the same papers, we shall identify elements
of the highest weight crystal $\B(\la)$ with semi-standard tableaux.

For later use, we recall the Kashiwara operator action on these
tableaux. We first read the boxes in the tableau through the \emph{far
eastern reading} and write down the boxes in \emph{tensor product
form.} That is, we read through each column from top to bottom starting
 from the rightmost column, continuing to the left. The following
diagram gives an example.
\begin{equation*}
\raisebox{-0.4\height}{
\begin{texdraw}%
\fontsize{6}{6}\selectfont
\textref h:C v:C
\drawdim em
\setunitscale 1.4
\move(0 3)\rlvec(5 0)
\move(0 2)\rlvec(5 0)\rlvec(0 1)
\move(0 1)\rlvec(4 0)\rlvec(0 2)
\move(0 0)\rlvec(2 0)\rlvec(0 2)
\move(0 0)\rlvec(0 3)
\move(1 0)\rlvec(0 3)
\move(3 1)\rlvec(0 2)
\move(2 1)\rlvec(0 2)
\htext(0.5 2.5){$1$}
\htext(1.5 2.5){$1$}
\htext(2.5 2.5){$1$}
\htext(3.5 2.5){$3$}
\htext(4.5 2.5){$\bar 1$}
\htext(0.5 1.5){$2$}
\htext(1.5 1.5){$3$}
\htext(2.5 1.5){$\bar 3$}
\htext(3.5 1.5){$0$}
\htext(0.5 0.5){$3$}
\htext(1.5 0.5){$0$}
\end{texdraw}%
}
=
\raisebox{-0.3\height}{
\begin{texdraw}%
\fontsize{6}{6}\selectfont
\textref h:C v:C
\drawdim em
\setunitscale 1.4
\move(0 0)\rlvec(1 0)\rlvec(0 1)\rlvec(-1 0)\rlvec(0 -1)
\htext(0.5 0.5){$\bar 1$}
\end{texdraw}%
}
\otimes
\raisebox{-0.3\height}{%
\begin{texdraw}%
\fontsize{6}{6}\selectfont
\textref h:C v:C
\drawdim em
\setunitscale 1.4
\move(0 0)\rlvec(1 0)\rlvec(0 1)\rlvec(-1 0)\rlvec(0 -1)
\htext(0.5 0.5){$3$}
\end{texdraw}%
}
\otimes
\raisebox{-0.3\height}{%
\begin{texdraw}%
\fontsize{6}{6}\selectfont
\textref h:C v:C
\drawdim em
\setunitscale 1.4
\move(0 0)\rlvec(1 0)\rlvec(0 1)\rlvec(-1 0)\rlvec(0 -1)
\htext(0.5 0.5){$0$}
\end{texdraw}%
}
\otimes
\raisebox{-0.3\height}{%
\begin{texdraw}%
\fontsize{6}{6}\selectfont
\textref h:C v:C
\drawdim em
\setunitscale 1.4
\move(0 0)\rlvec(1 0)\rlvec(0 1)\rlvec(-1 0)\rlvec(0 -1)
\htext(0.5 0.5){$1$}
\end{texdraw}%
}
\otimes
\raisebox{-0.3\height}{%
\begin{texdraw}%
\fontsize{6}{6}\selectfont
\textref h:C v:C
\drawdim em
\setunitscale 1.4
\move(0 0)\rlvec(1 0)\rlvec(0 1)\rlvec(-1 0)\rlvec(0 -1)
\htext(0.5 0.5){$\bar 3$}
\end{texdraw}%
}
\otimes
\raisebox{-0.3\height}{%
\begin{texdraw}%
\fontsize{6}{6}\selectfont
\textref h:C v:C
\drawdim em
\setunitscale 1.4
\move(0 0)\rlvec(1 0)\rlvec(0 1)\rlvec(-1 0)\rlvec(0 -1)
\htext(0.5 0.5){$1$}
\end{texdraw}%
}
\otimes
\raisebox{-0.3\height}{%
\begin{texdraw}%
\fontsize{6}{6}\selectfont
\textref h:C v:C
\drawdim em
\setunitscale 1.4
\move(0 0)\rlvec(1 0)\rlvec(0 1)\rlvec(-1 0)\rlvec(0 -1)
\htext(0.5 0.5){$3$}
\end{texdraw}%
}
\otimes
\raisebox{-0.3\height}{%
\begin{texdraw}%
\fontsize{6}{6}\selectfont
\textref h:C v:C
\drawdim em
\setunitscale 1.4
\move(0 0)\rlvec(1 0)\rlvec(0 1)\rlvec(-1 0)\rlvec(0 -1)
\htext(0.5 0.5){$0$}
\end{texdraw}%
}
\otimes
\raisebox{-0.3\height}{%
\begin{texdraw}%
\fontsize{6}{6}\selectfont
\textref h:C v:C
\drawdim em
\setunitscale 1.4
\move(0 0)\rlvec(1 0)\rlvec(0 1)\rlvec(-1 0)\rlvec(0 -1)
\htext(0.5 0.5){$1$}
\end{texdraw}%
}
\otimes
\raisebox{-0.3\height}{%
\begin{texdraw}%
\fontsize{6}{6}\selectfont
\textref h:C v:C
\drawdim em
\setunitscale 1.4
\move(0 0)\rlvec(1 0)\rlvec(0 1)\rlvec(-1 0)\rlvec(0 -1)
\htext(0.5 0.5){$2$}
\end{texdraw}%
}
\otimes
\raisebox{-0.3\height}{%
\begin{texdraw}%
\fontsize{6}{6}\selectfont
\textref h:C v:C
\drawdim em
\setunitscale 1.4
\move(0 0)\rlvec(1 0)\rlvec(0 1)\rlvec(-1 0)\rlvec(0 -1)
\htext(0.5 0.5){$3$}
\end{texdraw}%
}
\end{equation*}
Then, we apply the \emph{tensor product rule} to decide on which box to
apply $\fit$ or $\eit$ to. After application of Kashiwara operator to
one of the boxes, these are gathered back into the original form.

In practice, the tensor product rule on multiple tensors can be applied
through calculation of the \emph{$i$-signature}. This is done as
follows.
\begin{enumerate}
\item
First, under each tensor component $x$, write down $\veps_i(x)$-many 1s
followed by $\vphi_i(x)$-many 0s.
\item
Then, from the long sequence of mixed 0s and 1s, successively cancel
out every occurrence of (0,1) pair until we arrive at a sequence of 1s
followed by 0s, reading from left to right. This is called the
$i$-signature of the whole tensor product form.
\item
To apply $\fit$ to the whole product, apply it to the single tensor
component corresponding to the leftmost 0 remaining in the
$i$-signature. If no 0 remains, the result of $\fit$ action is set to
zero.
\item
Similarly, for $\eit$, apply it to the component corresponding to the
rightmost 1, or set it to zero when no 1 remains.
\end{enumerate}
We wish to restrict the set of dominant integral weights $P^+$ slightly
for some of the classical types.
\begin{itemize}
\item
$A_n$ case: $\rp := P^+$.
\item
$B_n$ case: $\rp := \{\la\in P^+\mid \text{$\la(h_n)$ is even}\}$.
\item
$C_n$ case: $\rp := P^+$.
\item
$D_{n+1}$ case: $\rp := \{\la\in P^+\mid \la(h_{n}) = \la(h_{n+1})\}$.
\item
$G_2$ case: $\rp := P^+$.
\end{itemize}
Notice that for $\la\in\rp$, elements of $\B(\la)$ become the most
\emph{generic} tableaux, in the sense that they do not involve any
half-size boxes or other complications. It is also clear that given any
$\la\in P^+$, we may always find a larger $\mu\in\rp$, that is, one
such that $\mu-\la\in P^+$.

We borrow the notion of \emph{large} semi-standard tableaux
 from~\cite{cliff}. For the remainder of this paper, the top row of a
tableau shall always be counted as the first row.
\begin{definition}
A semi-standard tableau $T$ of shape $\la\in\rp$ is \emph{large} if it
consists of $n$ non-empty rows, and if for each $1 \leq i \leq n$, the
number of $i$-boxes in the $i$-th row is strictly greater than the
number of all boxes in the $(i+1)$-th row. In particular the $n$-th row
of $T$ contains at least one $n$-box.

For each finite type, denote by $\calT(\la)^L$, the set of all large
semi-standard tableaux of shape $\la$.
\end{definition}
Once again, we remind readers that we are giving this definition for
each of the types $A_n$, $B_n$, $C_n$, $D_{n+1}$, and $G_2$. The $n$
appearing in the definition is meant to be the same $n$ used as
subscripts for the algebra types, with $n=2$ for the $G_2$ case.

In Figure~\ref{tbl:1}, for some of the finite types, we give examples
of semi-standard tableaux. The ones on the left column are large, and
the ones on the right are not large.
\begin{figure}
\centering
\begin{tabular}{rll}
$A_2$ case:
 &
 \raisebox{-0.4\height}{
 \begin{texdraw}%
 \drawdim in
 \arrowheadsize l:0.065 w:0.03
 \arrowheadtype t:F
 \fontsize{6}{6}\selectfont
 \textref h:C v:C
 \drawdim em
 \setunitscale 1.4
 \move(0 2)\rlvec(6 0)
 \move(0 1)\rlvec(6 0)\rlvec(0 1)
 \move(0 0)\rlvec(3 0)\rlvec(0 2)
 \move(0 0)\rlvec(0 2)
 \move(1 0)\rlvec(0 2)
 \move(2 0)\rlvec(0 2)
 \move(5 1)\rlvec(0 1)
 \move(4 1)\rlvec(0 1)
 \htext(0.5 1.5){$1$}
 \htext(1.5 1.5){$1$}
 \htext(2.5 1.5){$1$}
 \htext(3.5 1.5){$1$}
 \htext(4.5 1.5){$1$}
 \htext(5.5 1.5){$3$}
 \htext(0.5 0.5){$2$}
 \htext(1.5 0.5){$2$}
 \htext(2.5 0.5){$3$}
 \end{texdraw}%
 }
 &
 \raisebox{-0.4\height}{
 \begin{texdraw}%
 \drawdim in
 \arrowheadsize l:0.065 w:0.03
 \arrowheadtype t:F
 \fontsize{6}{6}\selectfont
 \textref h:C v:C
 \drawdim em
 \setunitscale 1.4
 \move(0 2)\rlvec(6 0)
 \move(0 1)\rlvec(6 0)\rlvec(0 1)
 \move(0 0)\rlvec(3 0)\rlvec(0 2)
 \move(0 0)\rlvec(0 2)
 \move(1 0)\rlvec(0 2)
 \move(2 0)\rlvec(0 2)
 \move(5 1)\rlvec(0 1)
 \move(4 1)\rlvec(0 1)
 \htext(0.5 1.5){$1$}
 \htext(1.5 1.5){$1$}
 \htext(2.5 1.5){$1$}
 \htext(3.5 1.5){$2$}
 \htext(4.5 1.5){$3$}
 \htext(5.5 1.5){$3$}
 \htext(0.5 0.5){$2$}
 \htext(1.5 0.5){$2$}
 \htext(2.5 0.5){$3$}
 \end{texdraw}%
 }\\[1.5em]
$B_3$ case:
 &
 \raisebox{-0.4\height}{
 \begin{texdraw}%
 \drawdim in
 \arrowheadsize l:0.065 w:0.03
 \arrowheadtype t:F
 \fontsize{6}{6}\selectfont
 \textref h:C v:C
 \drawdim em
 \setunitscale 1.4
 \move(0 3)\rlvec(8 0)
 \move(0 2)\rlvec(8 0)\rlvec(0 1)
 \move(0 1)\rlvec(4 0)\rlvec(0 2)
 \move(0 0)\rlvec(2 0)\rlvec(0 2)
 \move(0 0)\rlvec(0 3)
 \move(1 0)\rlvec(0 3)
 \move(3 1)\rlvec(0 2)
 \move(2 1)\rlvec(0 2)
 \move(5 2)\rlvec(0 1)
 \move(6 2)\rlvec(0 1)
 \move(7 2)\rlvec(0 1)
 \htext(0.5 2.5){$1$}
 \htext(1.5 2.5){$1$}
 \htext(2.5 2.5){$1$}
 \htext(3.5 2.5){$1$}
 \htext(4.5 2.5){$1$}
 \htext(5.5 2.5){$1$}
 \htext(6.5 2.5){$\bar 3$}
 \htext(7.5 2.5){$\bar 1$}
 \htext(0.5 1.5){$2$}
 \htext(1.5 1.5){$2$}
 \htext(2.5 1.5){$2$}
 \htext(3.5 1.5){$\bar 2$}
 \htext(0.5 0.5){$3$}
 \htext(1.5 0.5){$0$}
 \end{texdraw}%
 }
 &
 \raisebox{-0.4\height}{
 \begin{texdraw}%
 \drawdim in
 \arrowheadsize l:0.065 w:0.03
 \arrowheadtype t:F
 \fontsize{6}{6}\selectfont
 \textref h:C v:C
 \drawdim em
 \setunitscale 1.4
 \move(0 3)\rlvec(5 0)
 \move(0 2)\rlvec(5 0)\rlvec(0 1)
 \move(0 1)\rlvec(4 0)\rlvec(0 2)
 \move(0 0)\rlvec(2 0)\rlvec(0 2)
 \move(0 0)\rlvec(0 3)
 \move(1 0)\rlvec(0 3)
 \move(3 1)\rlvec(0 2)
 \move(2 1)\rlvec(0 2)
 \htext(0.5 2.5){$1$}
 \htext(1.5 2.5){$1$}
 \htext(2.5 2.5){$1$}
 \htext(3.5 2.5){$3$}
 \htext(4.5 2.5){$\bar 1$}
 \htext(0.5 1.5){$2$}
 \htext(1.5 1.5){$3$}
 \htext(2.5 1.5){$\bar 3$}
 \htext(3.5 1.5){$0$}
 \htext(0.5 0.5){$3$}
 \htext(1.5 0.5){$0$}
 \end{texdraw}%
 }\\[1.5em]
$C_4$ case:
 &
 \raisebox{-0.4\height}{
 \begin{texdraw}%
 \drawdim in
 \arrowheadsize l:0.065 w:0.03
 \arrowheadtype t:F
 \fontsize{6}{6}\selectfont
 \textref h:C v:C
 \drawdim em
 \setunitscale 1.4
 \move(0 4)\rlvec(10 0)
 \move(0 3)\rlvec(10 0)\rlvec(0 1)
 \move(0 2)\rlvec(7 0)\rlvec(0 2)
 \move(0 1)\rlvec(3 0)\rlvec(0 3)
 \move(0 0)\rlvec(2 0)\rlvec(0 4)
 \move(0 0)\rlvec(0 4)
 \move(1 0)\rlvec(0 4)
 \move(4 2)\rlvec(0 2)
 \move(5 2)\rlvec(0 2)
 \move(6 2)\rlvec(0 2)
 \move(8 3)\rlvec(0 1)
 \move(9 3)\rlvec(0 1)
 \htext(0.5 3.5){$1$}
 \htext(1.5 3.5){$1$}
 \htext(2.5 3.5){$1$}
 \htext(3.5 3.5){$1$}
 \htext(4.5 3.5){$1$}
 \htext(5.5 3.5){$1$}
 \htext(6.5 3.5){$1$}
 \htext(7.5 3.5){$1$}
 \htext(8.5 3.5){$4$}
 \htext(9.5 3.5){$\bar 1$}
 \htext(0.5 2.5){$2$}
 \htext(1.5 2.5){$2$}
 \htext(2.5 2.5){$2$}
 \htext(3.5 2.5){$2$}
 \htext(4.5 2.5){$2$}
 \htext(5.5 2.5){$\bar 3$}
 \htext(6.5 2.5){$\bar 2$}
 \htext(0.5 1.5){$3$}
 \htext(1.5 1.5){$3$}
 \htext(2.5 1.5){$3$}
 \htext(0.5 0.5){$4$}
 \htext(1.5 0.5){$\bar 4$}
 \end{texdraw}%
 }
 &
 \raisebox{-0.4\height}{
 \begin{texdraw}%
 \drawdim in
 \arrowheadsize l:0.065 w:0.03
 \arrowheadtype t:F
 \fontsize{6}{6}\selectfont
 \textref h:C v:C
 \drawdim em
 \setunitscale 1.4
 \move(0 4)\rlvec(7 0)
 \move(0 3)\rlvec(7 0)\rlvec(0 1)
 \move(0 2)\rlvec(3 0)\rlvec(0 2)
 \move(0 1)\rlvec(3 0)\rlvec(0 3)
 \move(0 0)\rlvec(2 0)\rlvec(0 4)
 \move(0 0)\rlvec(0 4)
 \move(1 0)\rlvec(0 4)
 \move(4 3)\rlvec(0 1)
 \move(5 3)\rlvec(0 1)
 \move(6 3)\rlvec(0 1)
 \htext(0.5 3.5){$1$}
 \htext(1.5 3.5){$1$}
 \htext(2.5 3.5){$1$}
 \htext(3.5 3.5){$1$}
 \htext(4.5 3.5){$2$}
 \htext(5.5 3.5){$4$}
 \htext(6.5 3.5){$\bar 3$}
 \htext(0.5 2.5){$2$}
 \htext(1.5 2.5){$2$}
 \htext(2.5 2.5){$3$}
 \htext(0.5 1.5){$3$}
 \htext(1.5 1.5){$3$}
 \htext(2.5 1.5){$\bar 4$}
 \htext(0.5 0.5){$\bar 4$}
 \htext(1.5 0.5){$\bar 4$}
 \end{texdraw}%
 }\\[2em]
$D_4$ case:
 &
 \raisebox{-0.4\height}{
 \begin{texdraw}%
 \drawdim in
 \arrowheadsize l:0.065 w:0.03
 \arrowheadtype t:F
 \fontsize{6}{6}\selectfont
 \textref h:C v:C
 \drawdim em
 \setunitscale 1.4
 \move(0 3)\rlvec(7 0)
 \move(0 2)\rlvec(7 0)\rlvec(0 1)
 \move(0 1)\rlvec(5 0)\rlvec(0 2)
 \move(0 0)\rlvec(2 0)\rlvec(0 3)
 \move(0 0)\rlvec(0 3)
 \move(1 0)\rlvec(0 3)
 \move(3 1)\rlvec(0 2)
 \move(4 1)\rlvec(0 2)
 \move(6 2)\rlvec(0 1)
 \htext(0.5 2.5){$1$}
 \htext(1.5 2.5){$1$}
 \htext(2.5 2.5){$1$}
 \htext(3.5 2.5){$1$}
 \htext(4.5 2.5){$1$}
 \htext(5.5 2.5){$1$}
 \htext(6.5 2.5){$4$}
 \htext(0.5 1.5){$2$}
 \htext(1.5 1.5){$2$}
 \htext(2.5 1.5){$2$}
 \htext(3.5 1.5){$\bar 4$}
 \htext(4.5 1.5){$\bar 4$}
 \htext(0.5 0.5){$3$}
 \htext(1.5 0.5){$\bar 3$}
 \end{texdraw}%
 }
 &
 \raisebox{-0.4\height}{
 \begin{texdraw}%
 \drawdim in
 \arrowheadsize l:0.065 w:0.03
 \arrowheadtype t:F
 \fontsize{6}{6}\selectfont
 \textref h:C v:C
 \drawdim em
 \setunitscale 1.4
 \move(0 3)\rlvec(6 0)
 \move(0 2)\rlvec(6 0)\rlvec(0 1)
 \move(0 1)\rlvec(5 0)\rlvec(0 2)
 \move(0 0)\rlvec(1 0)\rlvec(0 3)
 \move(0 0)\rlvec(0 3)
 \move(2 1)\rlvec(0 2)
 \move(3 1)\rlvec(0 2)
 \move(4 1)\rlvec(0 2)
 \htext(0.5 2.5){$1$}
 \htext(1.5 2.5){$1$}
 \htext(2.5 2.5){$2$}
 \htext(3.5 2.5){$2$}
 \htext(4.5 2.5){$3$}
 \htext(5.5 2.5){$\bar 4$}
 \htext(0.5 1.5){$3$}
 \htext(1.5 1.5){$4$}
 \htext(2.5 1.5){$4$}
 \htext(3.5 1.5){$\bar 3$}
 \htext(4.5 1.5){$\bar 2$}
 \htext(0.5 0.5){$\bar 4$}
 \end{texdraw}%
 }\\[1.5em]
$G_2$ case:
 &
 \raisebox{-0.4\height}{
 \begin{texdraw}%
 \drawdim in
 \arrowheadsize l:0.065 w:0.03
 \arrowheadtype t:F
 \fontsize{6}{6}\selectfont
 \textref h:C v:C
 \drawdim em
 \setunitscale 1.4
 \move(0 2)\rlvec(10 0)
 \move(0 1)\rlvec(10 0)\rlvec(0 1)
 \move(0 0)\rlvec(3 0)\rlvec(0 2)
 \move(0 0)\rlvec(0 2)
 \move(1 0)\rlvec(0 2)
 \move(2 0)\rlvec(0 2)
 \move(6 1)\rlvec(0 1)
 \move(7 1)\rlvec(0 1)
 \move(8 1)\rlvec(0 1)
 \move(9 1)\rlvec(0 1)
 \move(5 1)\rlvec(0 1)
 \move(4 1)\rlvec(0 1)
 \htext(0.5 1.5){$1$}
 \htext(1.5 1.5){$1$}
 \htext(2.5 1.5){$1$}
 \htext(3.5 1.5){$1$}
 \htext(9.5 1.5){$\bar 1$}
 \htext(8.5 1.5){$\bar 2$}
 \htext(7.5 1.5){$\bar 3$}
 \htext(6.5 1.5){$0$}
 \htext(5.5 1.5){$3$}
 \htext(4.5 1.5){$2$}
 \htext(0.5 0.5){$2$}
 \htext(1.5 0.5){$2$}
 \htext(2.5 0.5){$3$}
 \htext(1.5 0.5){$2$}
 \htext(2.5 0.5){$3$}
 \end{texdraw}%
 }
 &
 \raisebox{-0.4\height}{
 \begin{texdraw}%
 \drawdim in
 \arrowheadsize l:0.065 w:0.03
 \arrowheadtype t:F
 \fontsize{6}{6}\selectfont
 \textref h:C v:C
 \drawdim em
 \setunitscale 1.4
 \move(0 2)\rlvec(6 0)
 \move(0 1)\rlvec(6 0)\rlvec(0 1)
 \move(0 0)\rlvec(3 0)\rlvec(0 2)
 \move(0 0)\rlvec(0 2)
 \move(1 0)\rlvec(0 2)
 \move(2 0)\rlvec(0 2)
 \move(5 1)\rlvec(0 1)
 \move(4 1)\rlvec(0 1)
 \htext(0.5 1.5){$1$}
 \htext(1.5 1.5){$2$}
 \htext(2.5 1.5){$0$}
 \htext(3.5 1.5){$\bar{3}$}
 \htext(4.5 1.5){$\bar{3}$}
 \htext(5.5 1.5){$\bar{2}$}
 \htext(0.5 0.5){$2$}
 \htext(1.5 0.5){$3$}
 \htext(2.5 0.5){$\bar{2}$}
 \end{texdraw}%
 }
\end{tabular}\\[1em]
\caption{Large (left) and non-large (right) tableaux}\label{tbl:1}
\end{figure}

\section{The new crystal $\Tinf$}

Let us collect all large tableaux into one set (separately for each
finite type).
\begin{equation}
\calT^L = \cup_{\la\in\rp}\calT(\la)^L.
\end{equation}
We shall define an equivalence relation on this set.
\begin{definition}
Two tableaux $T_1, T_2 \in\calT^L$ are \emph{related} if for each
$1\leq i \leq n$ and $j \in J$ such that $j \succ i$, the number of
$j$-boxes appearing in the $i$-th rows of $T_1$ and $T_2$ are equal.
\end{definition}
It is trivial to verify that the above gives an equivalence relation.
We fix a notation
\begin{equation}
\Tinf := \TL/\sim
\end{equation}
for the set of equivalence classes. This section is devoted to
providing $\Tinf$ with a crystal structure.

Let us start with the Kashiwara operators.
\begin{lemma}\label{lem:3}
Fix an $i\in I$.
\begin{enumerate}
\item
If tableau $T$ is large, then $\fit T$ is never zero.
\item
Given any element of $\Tinf$, it is always possible to choose its
representative $T\in\TL$ in such a way that $\fit T$ is large.
\item\label{lem:3.3}
If $T_1, T_2\in\TL$ belong to the same equivalence class and $\fit T_1$
and $\fit T_2$ are both large, then $\fit T_1$ and $\fit T_2$ belong to
the same equivalence class.
\item
If tableau $T$ is large, then $\eit T$ is either zero or large.
\item
If $T_1, T_2\in\TL$ belong to the same equivalence class, then either
$\eit T_1$ and $\eit T_2$ are both zero, or $\eit T_1$ and $\eit T_2$
belong to the same equivalence class.
\end{enumerate}
\end{lemma}
\begin{proof}
1) It suffices to show that, after all cancelling out, at least one $0$
remains in the $i$-signature for $T$. Consider the rightmost $i$-block
on the $i$-th row of $T$. The signature to be written under it in the
tensor form of $T$ is $0$. Notice that the condition \emph{large}
guarantees it to be the lowest block in its column. Such careful
consideration of both the conditions \emph{large} and
\emph{semi-standard} for each of the finite types will show that the
signature $0$ under that block will not be cancelled out by signatures
 from blocks contained in any of the columns sitting to its left.

\smallskip\noindent 2)
Given any $T\in\TL$ such that $\bar{T} = b$, let us create a
\emph{larger} representative of $b$. First, construct a column
consisting of $i$ boxes, with $k$-box sitting on the $k$-th row ($1\leq
k \leq i$). Consider the rightmost $i$-box sitting on the $i$-th row of
$T$ and insert the constructed column to its left. It is clear that
this new tableau $T'$ is a (large) representative for $b$.

Now, during the proof of item~1 of this lemma, we saw that if we apply
$\fit$ to $T'$, it will act on either the rightmost $i$-block on the
$i$-th row of $T'$, or on one of the boxes sitting in columns to its
right. Due to the column we have inserted, neither case will affect the
largeness of $T'$, and hence the result is obtained.

\smallskip\noindent 3)
The tableaux $T_1$ and $T_2$ (or any other large tableaux) will take
the following form.
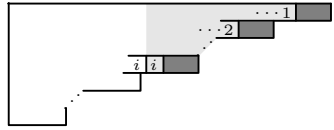
\begin{center}
\begin{texdraw}
\fontsize{6}{6}\selectfont
\textref h:C v:C
\drawdim em
\setunitscale 0.9
\move(0 0)\rlvec(-2 0)\rlvec(0 -1)\rlvec(2 0)\ifill f:0.5
\move(-3.3 -1)\rlvec(-2 0)\rlvec(0 -1)\rlvec(2 0)\ifill f:0.5
\move(-7.6 -3)\rlvec(-2 0)\rlvec(0 -1)\rlvec(2 0)\ifill f:0.5
\move(-10.6 0)\rlvec(0 -4)\rlvec(1 0)\rlvec(0 1)
\rlvec(2 0)\rlvec(1 1)\rlvec(1.3 0)\rlvec(0 1)
\rlvec(3.3 0)\rlvec(0 1)\ifill f:0.9
\move(0 0)\rlvec(-18.5 0)
\move(0 0)\rlvec(0 -1)\rlvec(-6.3 0)
\move(-2 0)\rlvec(0 -1)
\move(-3.3 -1)\rlvec(0 -1)\rlvec(-3.3 0)
\move(-5.3 -1)\rlvec(0 -1)
\move(-7.6 -3)\rlvec(-4 0)
\move(-7.6 -3)\rlvec(0 -1)\rlvec(-4.3 0)
\move(-9.6 -3)\rlvec(0 -1)
\move(-10.6 -3)\rlvec(0 -1)
\move(-10.9 -4)\rlvec(0 -1)\rlvec(-3.3 0)
\move(-15.2 -6)\rlvec(0 -1)\rlvec(-3.3 0)\rlvec(0 7)
\htext(-2.5 -0.5){$1$}
\htext(-3.5 -0.5){$\cdots$}
\htext(-5.8 -1.5){$2$}
\htext(-6.8 -1.5){$\cdots$}
\rtext td:45 (-7.1 -2.5){$\cdots$}
\rtext td:45 (-14.7 -5.5){$\cdots$}
\htext(-10.1 -3.5){$i$}
\htext(-11.1 -3.5){$i$}
\end{texdraw}
\end{center}
We already know from item~1 of this lemma that $\fit$ will not act on
any of the boxes contained in the unshaded part. It will act on either
the light shaded $i$-box or on the dark shaded part. Notice that, for
all finite types, $j$-boxes with $j < i$ do not contribute to
$i$-signatures, hence except for the $i$-box, none of the light shaded
part affects the $i$-signature for the two tableaux. Also, by
definition of the equivalence relation, the dark shaded parts will be
identical for the two. Thus all ingredients that become involved in the
$i$-signatures for the two tableaux are identical and hence $\fit$ will
act on two \emph{corresponding} boxes contained in the two tableaux.
This will result in the two tableaux being related even after $\fit$
action.

\smallskip\noindent 4) and 5)
These should be trivial to prove once proofs for previous items are
understood.
\end{proof}

It is now clear that, given $b\in\Tinf$ and $i\in I$, we may define
\begin{align}
\fit b &= \overline{\fit T} \in \Tinf,\label{eq:03}\\
\eit b &= \overline{\eit T} \in \Tinf \cup \{0\}
\end{align}
by choosing an appropriate representative $T$ for $b$.

\begin{lemma}
If $T_1\in\calT(\la_1)^L$ and $T_2\in\calT(\la_2)^L$ are related to
each other with $\wt(T_1) = \la_1 -\xi_1$ and $\wt(T_2) = \la_2
-\xi_2$, then $\xi_1 = \xi_2$.
\end{lemma}
Based on this trivial lemma, for $\bar{T}\in \Tinf$ with
$T\in\calT(\la)^L$, we can define
\begin{equation}
\wt(\bar{T}) = \wt(T) - \la.
\end{equation}
To complete the description of the crystal structure, it only remains
to define
\begin{align}
\veps_i(\bar{T}) &= \veps_i(T),\\
\vphi_i(\bar{T}) &= \veps_i(\bar{T}) + \wt(\bar{T})(h_i).\label{eq:07}
\end{align}
As before, these may be shown to be well-defined.

\begin{theorem}
The operator given by equations~\eqref{eq:03} to~\eqref{eq:07}, define
a crystal structure on $\Tinf$.
\end{theorem}
This theorem may be proved by a step by step checking of the definition
for an abstract crystal.

\section{Crystal isomorphism}\label{sec:4}

In this section, an isomorphism between crystal $\Binf$ and the crystal
$\Tinf$, constructed in the previous section, will be given. We start
by recalling the following theorem from~\cite{K91}.
\begin{theorem}\label{thm:6}
Let $\pi_{\la}: U_q^-(\mathfrak{g}) \rightarrow V(\la)$ be the
$U_q^-(\mathfrak{g})$-linear homomorphism sending 1 to $u_{\la}$. Then
\begin{enumerate}
\item
$\pi_{\la}(L(\infty)) = L(\la)$. Hence $\pi_{\la}$ induces the
surjective homomorphism
\begin{equation*}
\bpl:L(\infty)/qL(\infty) \rightarrow L(\la)/qL(\la).
\end{equation*}
\item
By $\bpl$, $\{b\in \Binf; \bpl(b)\neq 0\}$ is isomorphic to $\B(\la)$.
\item
$\fit\circ\bpl = \bpl\circ\fit$.

\item
If $b\in \Binf$ satisfies $\bpl(b)\neq 0$, then $\eit\bpl(b) =
\bpl(\eit b)$.
\end{enumerate}
\end{theorem}
We shall adopt the notation $\bpl$ introduced in this theorem. Let us
prepare for the definition of an explicit mapping from $\Binf$ to
$\Tinf$.
\begin{lemma}\hfill
\begin{enumerate}
\item
Given any $b\in\Binf$, there exists a $\la\in\rp$, for which $\bpl(b)$
is large.
\item
Given any $b\in\Binf$, if both $\bpl(b)$ and $\bar{\pi}_{\la'}(b)$ are
large, then the two belong to the same equivalence class. In
particular, any large highest weight elements $u_\la$ and $u_{\la'}$
are related.
\end{enumerate}
\end{lemma}
\begin{proof}
1) Simply put, depending on the distance of $b$ (in terms of $\fit$)
 from the highest element $u_\infty$, we can always choose $\la$ large
enough so that $\bpl(b)$ is still large. A more careful proof can be
written by modifying the proof of~\cite[Lemma~3.2]{cliff}.

\smallskip\noindent 2)
That any two large $u_\la$ and $u_{\la'}$ are related follows from the
definition for the equivalence relation. Starting from this point, we
may use induction with the help of Lemma~\ref{lem:3}~(\ref{lem:3.3}) to
obtain the result.
\end{proof}
We are now ready to define the mapping
\begin{equation}\label{eq:08}
\iso: \Binf \rightarrow \Tinf.
\end{equation}
Given any $b\in\Binf$, choose $\la\in\rp$ for which $\bpl(b)$ is large,
and set
\begin{equation}
\iso(b) = \overline{\bpl(b)}.
\end{equation}
The above lemma shows that this is well-defined. We thus arrive at one
of our main results.

\begin{theorem}
The mapping~\eqref{eq:08} is an isomorphism between $\Binf$ and
$\Tinf$.
\end{theorem}
\begin{proof}
It is easy to check that $\iso$ preserves $\veps_i$, $\vphi_i$, and
$\wt$. Notice that in the definition of $\iso$, the case $\bpl(b) = 0$
is never encountered. Hence items~3 and~4 of Theorem~\ref{thm:6} show
that $\iso$ is a strict crystal morphism. Also, since we know
$|\Binf_{-\xi}| = |\calT(\la)_{\la-\xi}|$ for all sufficiently large
$\la$ (Cor.~4.4.5 of~\cite{K91}), the mapping is bijective.
\end{proof}

\section{An explicit description of $\Binf$}\label{sec:5}

To achieve our final goal of giving an explicit description of $\Binf$
in terms of tableaux, it suffices to describe an explicit set of
representatives for $\Tinf = \TL/\sim$ and translate the various
operators on $\Tinf$ to that on the representative set.

\begin{definition}
A tableau $T\in\TL$ is \emph{marginally} large, if for $1\leq i \leq
n$, the number of $i$-boxes in the $i$-th row of $T$ is greater than
the number of all boxes in the $(i+1)$-th row by exactly one. In
particular, the $n$-th row of $T$ should contain one $n$-box.
\end{definition}

It is clear that the set of marginally large tableaux forms a set of
representatives for $\Tinf \cong \Binf$. In passing, we remark that the
difference of numbers considered above does not have to be \emph{one}
to obtain a representative set. It suffices to fix it to some
non-negative number for each row.

We shall now describe this representative set more explicitly for each
finite type. For each case, we shall present a set of alphabets to be
used inside the boxes forming the tableaux, together with an ordering
on the set. Next, a set of conditions that should be satisfied by the
tableaux is presented. The set of all tableaux subject to the given
conditions will be the set of representatives for $\Tinf \cong \Binf$.

These descriptions were obtained by considering all conditions defining
\emph{semi-standard} tableaux together with the condition
\emph{marginally large}. Trickiest part of the notion semi-standard
involves something called \emph{configuration}, but the condition large
ensures that no such configuration can occur, and we obtain a vast
simplification. The final description we give below are thus much
simpler than the definition for semi-standard tableaux.

After giving out the explicit representative sets, we shall describe
action of Kashiwara operators on these sets, in a manner which is
applicable commonly to all cases. We leave translation of other
operators $\wt$, $\veps_i$, and $\vphi_i$ to the readers.

Subsections~\ref{sec:51} to~\ref{sec:56} can be seen as the main
contribution of this paper.

\subsection{$A_n$ case}\label{sec:51}

Alphabet:
\begin{equation*}
J = \{ 1 \prec 2 \prec \cdots \prec n \prec n+1 \}.
\end{equation*}
Conditions:
\begin{enumerate}
\item
Tableau consists of $n$ rows.

\item
For $1\leq i \leq n$, the $i$-th row of the leftmost column is an
$i$-box.

\item
Box indices weakly increase (w.r.t. $\prec$) as we go to the right.

\item
For $1 \leq i \leq n$, the number of $i$-boxes in the $i$-th row is
larger than the total number of boxes appearing in the $(i+1)$-th row
by exactly one.
\end{enumerate}
\begin{example}
The set of representatives for $\Tinf$, in the $A_2$ case, consists of
all tableaux of the following form. The unshaded part must exist,
whereas the shaded part is optional with variable size.
\begin{equation*}
T = \raisebox{-0.5\height}{\ %
\begin{texdraw}
 \fontsize{6}{6}\selectfont
 \textref h:C v:C
 \drawdim em
 \setunitscale 1.35
 \move(6 0)\lvec(1 0)\lvec(1 -1)\lvec(6 -1)\lvec(6 0)\ifill f:0.8
 \move(0 0)\lvec(-2.5 0)\lvec(-2.5 -2)\lvec(0 -2)\lvec(0 0)\ifill f:0.8
 \move(6 0)\lvec(-3.5 0)
 \move(6 -1)\lvec(-3.5 -1)
 \move(0 -2)\lvec(-3.5 -2)
 \move(-3.5 0)\lvec(-3.5 -2)
 \htext(-3 -0.5){$1$}
 \htext(-0.55 -0.5){$1$}
 \htext(-1.7 -0.5){$\cdots$}

 \move(0 0)\rlvec(0 -2)
 \htext(0.5 -0.5){$1$}
 \move(1 0)\rlvec(0 -1)
 \htext(2.25 -0.5){$2\!\cdots\!2$}
 \move(3.5 0)\rlvec(0 -1)
 \htext(4.75 -0.5){$3\!\cdots\!3$}
 \move(6 0)\rlvec(0 -1)
 \htext(-1.25 -1.5){$3\!\cdots\!3$}
 \move(-2.5 -0.9)\rlvec(0 -1.1)
 \htext(-3 -1.5){$2$}
\end{texdraw}}
\end{equation*}
The element corresponding to highest weight element $u_\infty$ is
\begin{equation*}
T_{\infty} = \raisebox{-0.5\height}{\ %
\begin{texdraw}
 \fontsize{6}{6}\selectfont
 \textref h:C v:C
 \drawdim em
 \setunitscale 1.35
 \move(0 0)\lvec(2 0)\lvec(2 -1)\lvec(0 -1)
 \move(1 0)\lvec(1 -2)\lvec(0 -2)\lvec(0 0)
\htext(0.5 -0.5){$1$}
\htext(1.5 -0.5){$1$}
\htext(0.5 -1.5){$2$}
\end{texdraw}}\,.
\end{equation*}
\end{example}

\subsection{$B_n$ case}

Alphabet:
\begin{equation*}
J = \{ 1 \prec 2 \prec \cdots \prec n \prec 0 \prec \bar{n} \prec
\cdots \prec \bar{2} \prec \bar{1} \}.
\end{equation*}
Conditions:
\begin{enumerate}
\item
Tableau consists of $n$ rows.

\item
For $1\leq i \leq n$, the $i$-th row of the leftmost column is an
$i$-box.

\item
Box indices weakly increase (w.r.t. $\prec$) as we go to the right.

\item
For $1 \leq i \leq n$, the number of $i$-boxes in the $i$-th row is
larger than the total number of boxes appearing in the $(i+1)$-th row
by exactly one.

\item
All entries on the $i$-th row are less than or equal to $\bar{i}$
(w.r.t. $\prec$).

\item
Index $0$ appears at most once in each row.
\end{enumerate}

\begin{example}
The set of representatives for $\Tinf$, in the $B_3$ case, consists of
all tableaux of the following form. The unshaded part must exist,
whereas the shaded part is optional with variable size.
\begin{equation*}
T = \raisebox{-0.5\height}{\ %
\begin{texdraw}
 \fontsize{6}{6}\selectfont
 \textref h:C v:C
 \drawdim em
 \setunitscale 1.35
 \move(14.5 0)\lvec(14.5 -1)\lvec(1 -1)\lvec(1 0)\lvec(14.5 0)
 \ifill f:0.8
 \move(0 0)\lvec(0 -2)\lvec(-8.5 -2)\lvec(-8.5 0)\lvec(0 0)
 \ifill f:0.8
 \move(-9.5 0)\lvec(-9.5 -3)\lvec(-13 -3)\lvec(-13 0)\lvec(-9.5 0)
 \ifill f:0.8
 \move(14.5 0)\lvec(-14 0)
 \move(14.5 -1)\lvec(-14 -1)
 \move(0 -2)\lvec(-14 -2)
 \move(-9.5 -3)\lvec(-14 -3)
 \move(-14 0)\rlvec(0 -3)
 \htext(-12.75 -0.5){$1\!\cdots\!\phantom{1}$}
 \htext(-1.25 -0.5){$\phantom{1}\!\cdots\!1$}
 \move(0 0)\rlvec(0 -2)
 \htext(0.5 -0.5){$1$}
 \move(1 0)\rlvec(0 -1)
 \htext(2.25 -0.5){$2\!\cdots\!2$}
 \move(3.5 0)\rlvec(0 -1)
 \htext(4.75 -0.5){$3\!\cdots\!3$}
 \move(6 0)\rlvec(0 -1)
 \htext(6.5 -0.5){$0$}
 \move(7 0)\rlvec(0 -1)
 \htext(8.25 -0.5){$\bar{3}\!\cdots\!\bar{3}$}
 \move(9.5 0)\rlvec(0 -1)
 \htext(10.75 -0.5){$\bar{2}\!\cdots\!\bar{2}$}
 \move(12 0)\rlvec(0 -1)
 \htext(13.25 -0.5){$\bar{1}\!\cdots\!\bar{1}$}
 \move(14.5 0)\rlvec(0 -1)
 \htext(-12.75 -1.5){$2\!\cdots\!\phantom{2}$}
 \htext(-1.25 -1.5){$\bar{2}\!\cdots\!\bar{2}$}
 \move(-2.5 -0.9)\rlvec(0 -1.1)
 \htext(-3.75 -1.5){$\bar{3}\!\cdots\!\bar{3}$}
 \move(-5 -0.9)\rlvec(0 -1.1)
 \htext(-5.5 -1.5){$0$}
 \move(-6 -0.9)\rlvec(0 -1.1)
 \htext(-7.25 -1.5){$3\!\cdots\!3$}
 \move(-8.5 -0.9)\rlvec(0 -1.1)
 \htext(-9 -1.5){$2$}
 \move(-9.5 -0.9)\rlvec(0 -2.1)
 \htext(-10.75 -1.5){$\phantom{2}\!\cdots\!2$}
 \htext(-10.75 -2.5){$\bar{3}\!\cdots\!\bar{3}$}
 \move(-12 -1.9)\rlvec(0 -1.1)
 \htext(-12.5 -2.5){$0$}
 \move(-13 -1.9)\rlvec(0 -1.1)
 \htext(-13.5 -2.5){$3$}
\end{texdraw}}
\end{equation*}
The element corresponding to highest weight element $u_\infty$ is
\begin{equation*}
T_{\infty} = \raisebox{-0.5\height}{\ %
\begin{texdraw}
 \fontsize{6}{6}\selectfont
 \textref h:C v:C
 \drawdim em
 \setunitscale 1.35
 \move(1 0)\lvec(1 -3)\lvec(0 -3)\lvec(0 0)\lvec(3 0)
 \lvec(3 -1)\lvec(0 -1)
 \move(2 0)\lvec(2 -2)\lvec(0 -2)
\htext(0.5 -0.5){$1$}
\htext(1.5 -0.5){$1$}
\htext(2.5 -0.5){$1$}
\htext(0.5 -1.5){$2$}
\htext(1.5 -1.5){$2$}
\htext(0.5 -2.5){$3$}
\end{texdraw}}.
\end{equation*}
\end{example}

\subsection{$C_n$ case}

Alphabet:
\begin{equation*}
J = \{ 1 \prec 2 \prec \cdots \prec n \prec \bar{n} \prec \cdots \prec
\bar{2} \prec \bar{1} \}.
\end{equation*}
Conditions:
\begin{enumerate}
\item
Tableau consists of $n$ rows.

\item
For $1\leq i \leq n$, the $i$-th row of the leftmost column is an
$i$-box.

\item
Box indices weakly increase (w.r.t. $\prec$) as we go to the right.

\item
For $1 \leq i \leq n$, the number of $i$-boxes in the $i$-th row is
larger than the total number of boxes appearing in the $(i+1)$-th row
by exactly one.

\item
All entries on the $i$-th row are less than or equal to $\bar{i}$
(w.r.t. $\prec$).
\end{enumerate}
\begin{example}
The set of representatives for $\Tinf$, in the $C_3$ case, consists of
all tableaux of the following form. The unshaded part must exist,
whereas the shaded part is optional with variable size.
\begin{equation*}
T = \raisebox{-0.5\height}{\ %
\begin{texdraw}
 \fontsize{6}{6}\selectfont
 \textref h:C v:C
 \drawdim em
 \setunitscale 1.35
 \move(13.5 0)\lvec(1 0)\lvec(1 -1)\lvec(13.5 -1)\lvec(13.5 0)
 \ifill f:0.8
 \move(0 0)\lvec(-7.5 0)\lvec(-7.5 -2)\lvec(0 -2)\lvec(0 0)
 \ifill f:0.8
 \move(-8.5 0)\lvec(-8.5 -3)\lvec(-11 -3)\lvec(-11 0)\lvec(-8.5 0)
 \ifill f:0.8
 \move(13.5 0)\lvec(-12 0)
 \move(13.5 -1)\lvec(-12 -1)
 \move(0 -2)\lvec(-12 -2)
 \move(-9.5 -3)\lvec(-12 -3)
 \move(13.5 0)\rlvec(0 -1)
 \move(-8.5 -0.9)\rlvec(0 -2.1)\lvec(-12 -3)
 \move(-12 -3)\lvec(-12 0)
 \htext(-10.75 -0.5){$1\!\cdots\!\phantom{1}$}
 \htext(-1.25 -0.5){$\phantom{1}\!\cdots\!1$}
 \move(0 0)\rlvec(0 -2)
 \htext(0.5 -0.5){$1$}
 \move(1 0)\rlvec(0 -1)
 \htext(2.25 -0.5){$2\!\cdots\!2$}
 \move(3.5 0)\rlvec(0 -1)
 \htext(4.75 -0.5){$3\!\cdots\!3$}
 \move(6 0)\rlvec(0 -1)
 \htext(7.25 -0.5){$\bar{3}\!\cdots\!\bar{3}$}
 \move(8.5 0)\rlvec(0 -1)
 \htext(9.75 -0.5){$\bar{2}\!\cdots\!\bar{2}$}
 \move(11 0)\rlvec(0 -1)
 \htext(12.25 -0.5){$\bar{1}\!\cdots\!\bar{1}$}
 \htext(-1.25 -1.5){$\bar{2}\!\cdots\!\bar{2}$}
 \move(-2.5 -0.9)\rlvec(0 -1.1)
 \htext(-3.75 -1.5){$\bar{3}\!\cdots\!\bar{3}$}
 \move(-5 -0.9)\rlvec(0 -1.1)
 \htext(-6.25 -1.5){$3\!\cdots\!3$}
 \move(-7.5 -0.9)\rlvec(0 -1.1)
 \htext(-9 -1.5){$2$}
 \htext(-8 -1.5){$2$}
 \htext(-10.2 -1.5){$\cdots$}
 \move(-11 -3)\rlvec(0 1.1)
 \htext(-11.5 -1.5){$2$}
 \htext(-9.75 -2.5){$\bar{3}\!\cdots\!\bar{3}$}
 \htext(-11.5 -2.5){$3$}
\end{texdraw}}.
\end{equation*}
The element corresponding to highest weight element $u_\infty$ is
\begin{equation*}
T_{\infty} = \raisebox{-0.5\height}{\ %
\begin{texdraw}
 \fontsize{6}{6}\selectfont
 \textref h:C v:C
 \drawdim em
 \setunitscale 1.35
 \move(1 0)\lvec(1 -3)\lvec(0 -3)\lvec(0 0)\lvec(3 0)
 \lvec(3 -1)\lvec(0 -1)
 \move(2 0)\lvec(2 -2)\lvec(0 -2)
\htext(0.5 -0.5){$1$}
\htext(1.5 -0.5){$1$}
\htext(2.5 -0.5){$1$}
\htext(0.5 -1.5){$2$}
\htext(1.5 -1.5){$2$}
\htext(0.5 -2.5){$3$}
\end{texdraw}}.
\end{equation*}
\end{example}

\subsection{$D_{n+1}$ case}

Alphabet:
\begin{equation*}
J = \{ 1 \prec 2 \prec \cdots \prec n \prec
\genfrac{}{}{0pt}{}{\overline{n\!+\!1}}{n\!+\!1} \prec \bar{n} \prec
\cdots \prec \bar{2} \prec \bar{1} \}.
\end{equation*}
Conditions:
\begin{enumerate}
\item
Tableau consists of $n$ rows.

\item
For $1\leq i \leq n$, the $i$-th row of the leftmost column is an
$i$-box.

\item
Box indices weakly increase (w.r.t. $\prec$) as we go to the right.

\item
For $1 \leq i \leq n$, the number of $i$-boxes in the $i$-th row is
larger than the total number of boxes appearing in the $(i+1)$-th row
by exactly one.

\item
All entries on the $i$-th row are less than or equal to $\bar{i}$
(w.r.t. $\prec$).

\item
$n\!+\!1$ and $\overline{n\!+\!1}$ do not appear on the same row.
\end{enumerate}
\begin{example}
The set of representatives for $\Tinf$, in the $D_4$ case, consists of
all tableaux of the following form. The unshaded part must exist,
whereas the shaded part is optional with variable size. Either one of
$4$ or $\bar{4}$ may take the place of each of the letters $x$, $y$,
and $z$.
\begin{equation*}
T = \raisebox{-0.5\height}{\ %
\begin{texdraw}
 \fontsize{6}{6}\selectfont
 \textref h:C v:C
 \drawdim em
 \setunitscale 1.35
 \move(16 0)\lvec(1 0)\lvec(1 -1)\lvec(16 -1)\lvec(16 0)
 \ifill f:0.8
 \move(0 0)\lvec(-10 0)\lvec(-10 -2)\lvec(0 -2)\lvec(0 0)
 \ifill f:0.8
 \move(-16 0)\lvec(-11 0)\lvec(-11 -3)\lvec(-16 -3)\lvec(-16 0)
 \ifill f:0.8
 \move(16 0)\lvec(-17 0)
 \move(16 -1)\lvec(-17 -1)
 \move(0 -2)\lvec(-17 -2)
 \move(-17 0)\rlvec(0 -3)
 \move(16 0)\rlvec(0 -1)
 \htext(-15.75 -0.5){$1\!\cdots\!\phantom{1}$}
 \htext(-1.25 -0.5){$\phantom{1}\!\cdots\!1$}
 \move(0 0)\rlvec(0 -2)
 \htext(0.5 -0.5){$1$}
 \move(1 0)\rlvec(0 -1)
 \htext(2.25 -0.5){$2\!\cdots\!2$}
 \move(3.5 0)\rlvec(0 -1)
 \htext(4.75 -0.5){$3\!\cdots\!3$}
 \move(6 0)\rlvec(0 -1)
 \htext(7.25 -0.5){$x\!\cdots\!x$}

 \htext(9.75 -0.5){$\bar{3}\!\cdots\!\bar{3}$}
 \htext(12.25 -0.5){$\bar{2}\!\cdots\!\bar{2}$}
 \htext(14.75 -0.5){$\bar{1}\!\cdots\!\bar{1}$}
 \move(8.5 0)\rlvec(0 -1)
 \move(11 0)\rlvec(0 -1)
 \move(13.5 0)\rlvec(0 -1)
 \htext(-15.75 -1.5){$2\!\cdots\!\phantom{2}$}
 \htext(-1.25 -1.5){$\bar{2}\!\cdots\!\bar{2}$}
 \move(-2.5 -0.9)\rlvec(0 -1.1)
 \htext(-3.75 -1.5){$\bar{3}\!\cdots\!\bar{3}$}
 \htext(-6.25 -1.5){$y\!\cdots\!y$}
 \htext(-8.75 -1.5){$3\!\cdots\!3$}
 \htext(-10.5 -1.5){$2$}
 \move(-5 -0.9)\rlvec(0 -1.1)
 \move(-7.5 -0.9)\rlvec(0 -1.1)
 \move(-10 -0.9)\rlvec(0 -1.1)
 \move(-11 -0.9)\rlvec(0 -2.1)\rlvec(-6 0)
 \htext(-12.25 -1.5){$\phantom{2}\!\cdots\!2$}
 \htext(-12.25 -2.5){$\bar{3}\!\cdots\!\bar{3}$}
 \move(-13.5 -1.9)\rlvec(0 -1.1)
 \htext(-14.75 -2.5){$z\!\cdots\!z$}
 \move(-16 -3)\rlvec(0 1.1)
 \htext(-16.5 -2.5){$3$}
\end{texdraw}}.
\end{equation*}
The element corresponding to highest weight element $u_\infty$ is
\begin{equation*}
T_{\infty} = \raisebox{-0.5\height}{\ %
\begin{texdraw}
 \fontsize{6}{6}\selectfont
 \textref h:C v:C
 \drawdim em
 \setunitscale 1.35
 \move(1 0)\lvec(1 -3)\lvec(0 -3)\lvec(0 0)\lvec(3 0)
 \lvec(3 -1)\lvec(0 -1)
 \move(2 0)\lvec(2 -2)\lvec(0 -2)
\htext(0.5 -0.5){$1$}
\htext(1.5 -0.5){$1$}
\htext(2.5 -0.5){$1$}
\htext(0.5 -1.5){$2$}
\htext(1.5 -1.5){$2$}
\htext(0.5 -2.5){$3$}
\end{texdraw}}.
\end{equation*}
\end{example}

\subsection{$G_2$ case}

Alphabet:
\begin{equation*}
J = \{ 1 \prec 2 \prec 3 \prec 0 \prec \bar{3} \prec \bar{2} \prec
\bar{1} \}.
\end{equation*}
Conditions:
\begin{enumerate}
\item
Tableau consists of $2$ rows.

\item
For $1\leq i \leq 2$, the $i$-th row of the leftmost column is an
$i$-box.

\item
Box indices weakly increase (w.r.t. $\prec$) as we go to the right.

\item
For $1 \leq i \leq 2$, the number of $i$-boxes in the $i$-th row is
larger than the total number of boxes appearing in the $(i+1)$-th row
by exactly one.

\item
Only $2$ and $3$ appear as indices on the second row.

\item
Index $0$ appears at most once on the first row.
\end{enumerate}
\begin{example}
The set of representatives for $\Tinf$, in the $G_2$ case, consists of
all tableaux of the following form. The unshaded part must exist,
whereas the shaded part is optional with variable size.
\begin{equation*}
T = \raisebox{-0.5\height}{\ %
\begin{texdraw}
 \fontsize{6}{6}\selectfont
 \textref h:C v:C
 \drawdim em
 \setunitscale 1.35
 \move(14.5 0)\lvec(1 0)\lvec(1 -1)\lvec(14.5 -1)\lvec(14.5 0)
 \ifill f:0.8
 \move(0 0)\lvec(-2.5 0)\lvec(-2.5 -2)\lvec(0 -2)\lvec(0 0)
 \ifill f:0.8
 \move(14.5 0)\lvec(-3.5 0)
 \move(14.5 -1)\lvec(-3.5 -1)
 \move(0 -2)\lvec(-3.5 -2)
 \move(-3.5 0)\lvec(-3.5 -2)
 \htext(-3 -0.5){$1$}
 \htext(-0.55 -0.5){$1$}
 \htext(-1.7 -0.5){$\cdots$}

 \move(0 0)\rlvec(0 -2)
 \htext(0.5 -0.5){$1$}
 \move(1 0)\rlvec(0 -1)
 \htext(2.25 -0.5){$2\!\cdots\!2$}
 \move(3.5 0)\rlvec(0 -1)
 \htext(4.75 -0.5){$3\!\cdots\!3$}
 \move(6 0)\rlvec(0 -1)
 \htext(6.5 -0.5){$0$}
 \move(7 0)\rlvec(0 -1)
 \htext(8.25 -0.5){$\bar{3}\!\cdots\!\bar{3}$}
 \move(9.5 0)\rlvec(0 -1)
 \htext(10.75 -0.5){$\bar{2}\!\cdots\!\bar{2}$}
 \move(12 0)\rlvec(0 -1)
 \htext(13.25 -0.5){$\bar{1}\!\cdots\!\bar{1}$}
 \move(14.5 0)\rlvec(0 -1)
 \htext(-1.25 -1.5){$3\!\cdots\!3$}
 \move(-2.5 -0.9)\rlvec(0 -1.1)
 \htext(-3 -1.5){$2$}
\end{texdraw}}.
\end{equation*}
The element corresponding to highest weight element $u_\infty$ is
\begin{equation*}
T_{\infty} = \raisebox{-0.5\height}{\ %
\begin{texdraw}
 \fontsize{6}{6}\selectfont
 \textref h:C v:C
 \drawdim em
 \setunitscale 1.35
 \move(0 0)\lvec(2 0)\lvec(2 -1)\lvec(0 -1)
 \move(1 0)\lvec(1 -2)\lvec(0 -2)\lvec(0 0)
\htext(0.5 -0.5){$1$}
\htext(1.5 -0.5){$1$}
\htext(0.5 -1.5){$2$}
\end{texdraw}}\,.
\end{equation*}
\end{example}

\subsection{Kashiwara operators}\label{sec:56}

To apply $\fit$ to one of the representatives, we go through the
following procedure.
\begin{enumerate}
\item
Apply $\fit$ to the tableau as usual. That is, write it in tensor
product form, apply tensor product rule, and assemble back into
original tableau form.

\item
If the result is a large tableau, we are done. It is automatically
marginally large.

\item
If the result is not large, the $\fit$ was applied to the rightmost
$i$-box on the $i$-th row. Insert one column consisting of $i$ rows to
the left of the box $\fit$ acted upon. The added column should have a
$k$-box at the $k$-th row for $1\leq k \leq i$.
\end{enumerate}
To apply $\eit$ to one of the representatives, we go through the
following procedure.
\begin{enumerate}
\item
Apply $\eit$ to the tableau as usual.

\item
If the result is zero or a marginally large tableau, we are done.

\item
Otherwise, the result is large but not marginally large. The $\eit$
operator has acted on the box sitting to the right of the rightmost
$i$-box on the $i$-th row. Remove the column containing the changed
box. It will be of $i$ rows and have a $k$-box at the $k$-th row for
$1\leq k \leq i$.
\end{enumerate}

\begin{example}
In the Figures~\ref{fig:2} and~\ref{fig:3}, we illustrate the top part
of crystal $\mathcal T(\infty)$ for finite types $B_3$ and $G_2$. The
dark shaded blocks are the ones $\fit$ has acted upon, and the light
shadings show columns inserted to preserve largeness.
\end{example}

\begin{figure}
\centering
\begin{texdraw}%
\drawdim in
\arrowheadsize l:0.065 w:0.03
\arrowheadtype t:F
\fontsize{5}{5}\selectfont
\textref h:C v:C
\drawdim em
\setunitscale 1.2
\move(-1 0)
\bsegment
\move(0 3)\rlvec(3 0)
\move(0 2)\rlvec(3 0)\rlvec(0 1)
\move(0 1)\rlvec(2 0)\rlvec(0 2)
\move(0 0)\rlvec(1 0)\rlvec(0 3)
\move(0 3)\rlvec(0 -3)
\htext(0.5 2.5){$1$}
\htext(1.5 2.5){$1$}
\htext(2.5 2.5){$1$}
\htext(0.5 1.5){$2$}
\htext(1.5 1.5){$2$}
\htext(0.5 0.5){$3$}
\esegment
\move(-8 -6)
\bsegment
\move(3 3)\rlvec(1 0)\rlvec(0 -1)\rlvec(-1 0)\ifill f:0.6
\move(2 3)\rlvec(1 0)\rlvec(0 -1)\rlvec(-1 0)\ifill f:0.9
\move(0 3)\rlvec(4 0)
\move(0 2)\rlvec(4 0)\rlvec(0 1)
\move(0 1)\rlvec(2 0)\rlvec(0 2)
\move(0 0)\rlvec(1 0)\rlvec(0 3)
\move(0 3)\rlvec(0 -3)
\move(3 3)\rlvec(0 -1)
\htext(0.5 2.5){$1$}
\htext(3.5 2.5){$2$}
\htext(1.5 2.5){$1$}
\htext(2.5 2.5){$1$}
\htext(0.5 1.5){$2$}
\htext(1.5 1.5){$2$}
\htext(0.5 0.5){$3$}
\esegment
\move(3 -6)
\bsegment
\move(-2 2)\rlvec(1 0)\rlvec(0 -1)\rlvec(-1 0)\ifill f:0.6
\move(-3 1)\rlvec(1 0)\rlvec(0 2)\rlvec(-1 0)\ifill f:0.9
\move(-4 3)\rlvec(4 0)
\move(-4 2)\rlvec(4 0)\rlvec(0 1)
\move(-4 1)\rlvec(3 0)\rlvec(0 2)
\move(-4 0)\rlvec(1 0)\rlvec(0 3)
\move(-4 3)\rlvec(0 -3)
\move(-2 3)\rlvec(0 -2)
\htext(-0.5 2.5){$1$}
\htext(-1.5 2.5){$1$}
\htext(-2.5 2.5){$1$}
\htext(-3.5 2.5){$1$}
\htext(-1.5 1.5){$3$}
\htext(-2.5 1.5){$2$}
\htext(-3.5 1.5){$2$}
\htext(-3.5 0.5){$3$}
\esegment
\move(10 -6)
\bsegment
\move(-3 1)\rlvec(1 0)\rlvec(0 -1)\rlvec(-1 0)\ifill f:0.6
\move(-4 0)\rlvec(1 0)\rlvec(0 3)\rlvec(-1 0)\ifill f:0.9
\move(-4 3)\rlvec(4 0)
\move(-4 2)\rlvec(4 0)\rlvec(0 1)
\move(-4 1)\rlvec(3 0)\rlvec(0 2)
\move(-4 0)\rlvec(2 0)\rlvec(0 3)
\move(-4 3)\rlvec(0 -3)
\move(-3 3)\rlvec(0 -3)
\htext(-0.5 2.5){$1$}
\htext(-1.5 2.5){$1$}
\htext(-2.5 2.5){$1$}
\htext(-3.5 2.5){$1$}
\htext(-1.5 1.5){$2$}
\htext(-2.5 1.5){$2$}
\htext(-3.5 1.5){$2$}
\htext(-2.5 0.5){$0$}
\htext(-3.5 0.5){$3$}
\esegment
\move(-15 -12)
\bsegment
\move(3 3)\rlvec(1 0)\rlvec(0 -1)\rlvec(-1 0)\ifill f:0.6
\move(2 3)\rlvec(1 0)\rlvec(0 -1)\rlvec(-1 0)\ifill f:0.9
\move(0 3)\rlvec(5 0)
\move(0 2)\rlvec(5 0)\rlvec(0 1)
\move(0 1)\rlvec(2 0)\rlvec(0 2)
\move(0 0)\rlvec(1 0)\rlvec(0 3)
\move(0 3)\rlvec(0 -3)
\move(3 3)\rlvec(0 -1)
\move(4 3)\rlvec(0 -1)
\htext(0.5 2.5){$1$}
\htext(1.5 2.5){$1$}
\htext(2.5 2.5){$1$}
\htext(3.5 2.5){$2$}
\htext(4.5 2.5){$2$}
\htext(0.5 1.5){$2$}
\htext(1.5 1.5){$2$}
\htext(0.5 0.5){$3$}
\esegment
\move(-8 -12)
\bsegment
\move(3 3)\rlvec(1 0)\rlvec(0 -1)\rlvec(-1 0)\ifill f:0.6
\move(0 3)\rlvec(4 0)
\move(0 2)\rlvec(4 0)\rlvec(0 1)
\move(0 1)\rlvec(2 0)\rlvec(0 2)
\move(0 0)\rlvec(1 0)\rlvec(0 3)
\move(0 3)\rlvec(0 -3)
\move(3 3)\rlvec(0 -1)
\htext(0.5 2.5){$1$}
\htext(1.5 2.5){$1$}
\htext(2.5 2.5){$1$}
\htext(3.5 2.5){$3$}
\htext(0.5 1.5){$2$}
\htext(1.5 1.5){$2$}
\htext(0.5 0.5){$3$}
\esegment
\move(3.5 -12)
\bsegment
\move(-4 1)\rlvec(1 0)\rlvec(0 -1)\rlvec(-1 0)\ifill f:0.6
\move(-5 0)\rlvec(1 0)\rlvec(0 3)\rlvec(-1 0)\ifill f:0.9
\move(-1 3)\rlvec(1 0)\rlvec(0 -1)\rlvec(-1 0)\ifill f:0.6
\move(-2 2)\rlvec(1 0)\rlvec(0 1)\rlvec(-1 0)\ifill f:0.9
\move(-5 3)\rlvec(5 0)
\move(-5 2)\rlvec(5 0)\rlvec(0 1)
\move(-5 1)\rlvec(3 0)\rlvec(0 2)
\move(-5 0)\rlvec(2 0)\rlvec(0 3)
\move(-4 3)\rlvec(0 -3)
\move(-5 3)\rlvec(0 -3)
\move(-1 3)\rlvec(0 -1)
\htext(-0.5 2.5){$2$}
\htext(-1.5 2.5){$1$}
\htext(-2.5 2.5){$1$}
\htext(-3.5 2.5){$1$}
\htext(-4.5 2.5){$1$}
\htext(-2.5 1.5){$2$}
\htext(-3.5 1.5){$2$}
\htext(-4.5 1.5){$2$}
\htext(-3.5 0.5){$0$}
\htext(-4.5 0.5){$3$}
\esegment
\move(6 -12)
\bsegment
\move(3 2)\rlvec(1 0)\rlvec(0 -1)\rlvec(-1 0)\ifill f:0.6
\move(2 1)\rlvec(1 0)\rlvec(0 2)\rlvec(-1 0)\ifill f:0.9
\move(0 3)\rlvec(5 0)
\move(0 2)\rlvec(5 0)\rlvec(0 1)
\move(0 1)\rlvec(4 0)\rlvec(0 2)
\move(0 0)\rlvec(2 0)\rlvec(0 3)
\move(0 3)\rlvec(0 -3)
\move(3 3)\rlvec(0 -2)
\move(1 3)\rlvec(0 -3)
\htext(0.5 2.5){$1$}
\htext(1.5 2.5){$1$}
\htext(2.5 2.5){$1$}
\htext(3.5 2.5){$1$}
\htext(4.5 2.5){$1$}
\htext(0.5 1.5){$2$}
\htext(1.5 1.5){$2$}
\htext(2.5 1.5){$2$}
\htext(3.5 1.5){$3$}
\htext(0.5 0.5){$3$}
\htext(1.5 0.5){$0$}
\esegment
\move(13 -12)
\bsegment
\move(1 1)\rlvec(1 0)\rlvec(0 -1)\rlvec(-1 0)\ifill f:0.6
\move(0 3)\rlvec(4 0)
\move(0 2)\rlvec(4 0)\rlvec(0 1)
\move(0 1)\rlvec(3 0)\rlvec(0 2)
\move(0 0)\rlvec(2 0)\rlvec(0 3)
\move(0 3)\rlvec(0 -3)
\move(1 3)\rlvec(0 -3)
\htext(0.5 2.5){$1$}
\htext(1.5 2.5){$1$}
\htext(2.5 2.5){$1$}
\htext(3.5 2.5){$1$}
\htext(0.5 1.5){$2$}
\htext(1.5 1.5){$2$}
\htext(2.5 1.5){$2$}
\htext(0.5 0.5){$3$}
\htext(1.5 0.5){$\bar 3$}
\esegment
\move(-4.5 -18)
\bsegment
\move(-1 3)\rlvec(1 0)\rlvec(0 -1)\rlvec(-1 0)\ifill f:0.6
\move(-2 3)\rlvec(1 0)\rlvec(0 -1)\rlvec(-1 0)\ifill f:0.9
\move(-5 3)\rlvec(5 0)
\move(-5 2)\rlvec(5 0)\rlvec(0 1)
\move(-5 1)\rlvec(3 0)\rlvec(0 2)
\move(-5 0)\rlvec(1 0)\rlvec(0 3)
\move(-5 3)\rlvec(0 -3)
\move(-3 3)\rlvec(0 -2)
\move(-1 3)\rlvec(0 -1)
\htext(-0.5 2.5){$2$}
\htext(-1.5 2.5){$1$}
\htext(-2.5 2.5){$1$}
\htext(-3.5 2.5){$1$}
\htext(-4.5 2.5){$1$}
\htext(-2.5 1.5){$3$}
\htext(-3.5 1.5){$2$}
\htext(-4.5 1.5){$2$}
\htext(-4.5 0.5){$3$}
\esegment
\move(3.5 -18)
\bsegment
\move(-3 2)\rlvec(1 0)\rlvec(0 -1)\rlvec(-1 0)\ifill f:0.6
\move(-4 1)\rlvec(1 0)\rlvec(0 2)\rlvec(-1 0)\ifill f:0.9
\move(-5 3)\rlvec(5 0)
\move(-5 2)\rlvec(5 0)\rlvec(0 1)
\move(-5 1)\rlvec(4 0)\rlvec(0 2)
\move(-5 0)\rlvec(1 0)\rlvec(0 3)
\move(-5 3)\rlvec(0 -3)
\move(-2 3)\rlvec(0 -2)
\move(-3 3)\rlvec(0 -2)
\htext(-0.5 2.5){$1$}
\htext(-1.5 2.5){$1$}
\htext(-2.5 2.5){$1$}
\htext(-3.5 2.5){$1$}
\htext(-4.5 2.5){$1$}
\htext(-1.5 1.5){$3$}
\htext(-2.5 1.5){$3$}
\htext(-3.5 1.5){$2$}
\htext(-4.5 1.5){$2$}
\htext(-4.5 0.5){$3$}
\esegment
\move(11 -18)
\bsegment
\move(-2 2)\rlvec(1 0)\rlvec(0 -1)\rlvec(-1 0)\ifill f:0.6
\move(-4 3)\rlvec(4 0)
\move(-4 2)\rlvec(4 0)\rlvec(0 1)
\move(-4 1)\rlvec(3 0)\rlvec(0 2)
\move(-4 0)\rlvec(1 0)\rlvec(0 3)
\move(-4 3)\rlvec(0 -3)
\move(-2 3)\rlvec(0 -2)
\htext(-0.5 2.5){$1$}
\htext(-1.5 2.5){$1$}
\htext(-2.5 2.5){$1$}
\htext(-3.5 2.5){$1$}
\htext(-1.5 1.5){$0$}
\htext(-2.5 1.5){$2$}
\htext(-3.5 1.5){$2$}
\htext(-3.5 0.5){$3$}
\esegment
\move(0 0)
\bsegment
\htext(-5.1 -1.3){$1$}
\htext(0.5 -1.3){$2$}
\htext(5.1 -1.3){$3$}
\esegment
\move(0 -3.7)
\bsegment
\htext(-9.8 -3.5){$1$}
\htext(-6.5 -3.5){$2$}
\htext(-3.5 -3.5){$3$}
\htext(-1.5 -3.5){$1$}
\htext(0.5 -3.5){$2$}
\htext(3.5 -3.5){$3$}
\htext(5.5 -3.5){$1$}
\htext(8.7 -3.5){$2$}
\htext(11.9 -3.5){$3$}
\esegment
\move(0 1)
\move(0 -14)
\move(-4.2 -1.1)\ravec(-1 -1)
\move(0 -1.1)\ravec(0 -1)
\move(4.2 -1.1)\ravec(1 -1)
\move(-9 -7.1)\ravec(-1 -1)
\move(-6 -7.1)\ravec(0 -1)
\move(-3 -7.1)\ravec(1 -1)
\move(5.2 -7.1)\ravec(-1 -1)
\move(8.2 -7.1)\ravec(0 -1)
\move(11.2 -7.1)\ravec(1 -1)
\move(3 -7.1)\ravec(3 -7)
\move(-1 -7.1)\ravec(-3 -7)
\htext(-14 -13){$\vdots$}
\htext(14 -13){$\vdots$}
\htext(-8 -19){$\vdots$}
\htext(0 -19){$\vdots$}
\htext(8 -19){$\vdots$}
\lpatt(0.05 0.15)
\move(0 -7.1)\rlvec(0 -1)
\move(0 -13.1)\ravec(0 -1)
\end{texdraw}
\caption{Crystal $\Tinf$ for type $B_3$}\label{fig:2}
\end{figure}

\begin{figure}
\centering
\begin{texdraw}%
\drawdim in
\arrowheadsize l:0.065 w:0.03
\arrowheadtype t:F
\fontsize{5}{5}\selectfont
\textref h:C v:C
\drawdim em
\setunitscale 1.2
\move(-1 0)
\bsegment
\move(0 2)\rlvec(2 0)
\move(0 1)\rlvec(2 0)\rlvec(0 1)
\move(0 0)\rlvec(1 0)\rlvec(0 2)
\move(0 2)\rlvec(0 -2)
\htext(0.5 1.5){$1$}
\htext(1.5 1.5){$1$}
\htext(0.5 0.5){$2$}
\esegment
\move(-8 -5)
\bsegment
\move(2 2)\rlvec(1 0)\rlvec(0 -1)\rlvec(-1 0)\ifill f:0.6
\move(1 2)\rlvec(1 0)\rlvec(0 -1)\rlvec(-1 0)\ifill f:0.9
\move(0 2)\rlvec(3 0)
\move(0 1)\rlvec(3 0)\rlvec(0 1)
\move(0 0)\rlvec(1 0)\rlvec(0 2)
\move(0 0)\rlvec(0 2)
\move(2 1)\rlvec(0 1)
\htext(0.5 1.5){$1$}
\htext(1.5 1.5){$1$}
\htext(2.5 1.5){$2$}
\htext(0.5 0.5){$2$}
\esegment
\move(5 -5)
\bsegment
\move(1 1)\rlvec(1 0)\rlvec(0 -1)\rlvec(-1 0)\ifill f:0.6
\move(0 0)\rlvec(1 0)\rlvec(0 2)\rlvec(-1 0)\ifill f:0.9
\move(0 2)\rlvec(3 0)
\move(0 1)\rlvec(3 0)\rlvec(0 1)
\move(0 0)\rlvec(2 0)\rlvec(0 2)
\move(1 0)\rlvec(0 2)
\move(0 0)\rlvec(0 2)
\htext(0.5 1.5){$1$}
\htext(1.5 1.5){$1$}
\htext(2.5 1.5){$1$}
\htext(0.5 0.5){$2$}
\htext(1.5 0.5){$3$}
\esegment
\move(-12 -10)
\bsegment
\move(2 2)\rlvec(1 0)\rlvec(0 -1)\rlvec(-1 0)\ifill f:0.6
\move(1 2)\rlvec(1 0)\rlvec(0 -1)\rlvec(-1 0)\ifill f:0.9
\move(0 2)\rlvec(4 0)
\move(0 1)\rlvec(4 0)\rlvec(0 1)
\move(0 0)\rlvec(1 0)\rlvec(0 2)
\move(0 0)\rlvec(0 2)
\move(2 1)\rlvec(0 1)
\move(3 1)\rlvec(0 1)
\htext(0.5 1.5){$1$}
\htext(1.5 1.5){$1$}
\htext(2.5 1.5){$2$}
\htext(3.5 1.5){$2$}
\htext(0.5 0.5){$2$}
\esegment
\move(-5 -10)
\bsegment
\move(2 2)\rlvec(1 0)\rlvec(0 -1)\rlvec(-1 0)\ifill f:0.6
\move(0 2)\rlvec(3 0)
\move(0 1)\rlvec(3 0)\rlvec(0 1)
\move(0 0)\rlvec(1 0)\rlvec(0 2)
\move(0 0)\rlvec(0 2)
\move(2 1)\rlvec(0 1)
\htext(0.5 1.5){$1$}
\htext(1.5 1.5){$1$}
\htext(2.5 1.5){$3$}
\htext(0.5 0.5){$2$}
\esegment
\move(1 -10)
\bsegment
\move(3 2)\rlvec(1 0)\rlvec(0 -1)\rlvec(-1 0)\ifill f:0.6
\move(2 2)\rlvec(1 0)\rlvec(0 -1)\rlvec(-1 0)\ifill f:0.9
\move(0 2)\rlvec(4 0)
\move(0 1)\rlvec(4 0)\rlvec(0 1)
\move(0 0)\rlvec(2 0)\rlvec(0 2)
\move(0 0)\rlvec(0 2)
\move(1 0)\rlvec(0 2)
\move(3 1)\rlvec(0 1)
\htext(0.5 1.5){$1$}
\htext(1.5 1.5){$1$}
\htext(2.5 1.5){$1$}
\htext(3.5 1.5){$2$}
\htext(0.5 0.5){$2$}
\htext(1.5 0.5){$3$}
\esegment
\move(8 -10)
\bsegment
\move(1 1)\rlvec(1 0)\rlvec(0 -1)\rlvec(-1 0)\ifill f:0.6
\move(0 0)\rlvec(1 0)\rlvec(0 2)\rlvec(-1 0)\ifill f:0.9
\move(0 2)\rlvec(4 0)
\move(0 1)\rlvec(4 0)\rlvec(0 1)
\move(0 0)\rlvec(3 0)\rlvec(0 2)
\move(2 0)\rlvec(0 2)
\move(1 0)\rlvec(0 2)
\move(0 0)\rlvec(0 2)
\htext(0.5 1.5){$1$}
\htext(1.5 1.5){$1$}
\htext(2.5 1.5){$1$}
\htext(3.5 1.5){$1$}
\htext(0.5 0.5){$2$}
\htext(1.5 0.5){$3$}
\htext(2.5 0.5){$3$}
\esegment
\move(-19.5 -15)
\bsegment
\move(2 2)\rlvec(1 0)\rlvec(0 -1)\rlvec(-1 0)\ifill f:0.6
\move(1 2)\rlvec(1 0)\rlvec(0 -1)\rlvec(-1 0)\ifill f:0.9
\move(0 2)\rlvec(5 0)
\move(0 1)\rlvec(5 0)\rlvec(0 1)
\move(0 0)\rlvec(1 0)\rlvec(0 2)
\move(0 0)\rlvec(0 2)
\move(2 1)\rlvec(0 1)
\move(3 1)\rlvec(0 1)
\move(4 1)\rlvec(0 1)
\htext(0.5 1.5){$1$}
\htext(1.5 1.5){$1$}
\htext(2.5 1.5){$2$}
\htext(3.5 1.5){$2$}
\htext(4.5 1.5){$2$}
\htext(0.5 0.5){$2$}
\esegment
\move(-13 -15)
\bsegment
\move(3 2)\rlvec(1 0)\rlvec(0 -1)\rlvec(-1 0)\ifill f:0.6
\move(0 2)\rlvec(4 0)
\move(0 1)\rlvec(4 0)\rlvec(0 1)
\move(0 0)\rlvec(1 0)\rlvec(0 2)
\move(0 0)\rlvec(0 2)
\move(2 1)\rlvec(0 1)
\move(3 1)\rlvec(0 1)
\htext(0.5 1.5){$1$}
\htext(1.5 1.5){$1$}
\htext(2.5 1.5){$2$}
\htext(3.5 1.5){$3$}
\htext(0.5 0.5){$2$}
\esegment
\move(-7.5 -15)
\bsegment
\move(2 2)\rlvec(1 0)\rlvec(0 -1)\rlvec(-1 0)\ifill f:0.6
\move(0 2)\rlvec(3 0)
\move(0 1)\rlvec(3 0)\rlvec(0 1)
\move(0 0)\rlvec(1 0)\rlvec(0 2)
\move(0 0)\rlvec(0 2)
\move(2 1)\rlvec(0 1)
\htext(0.5 1.5){$1$}
\htext(1.5 1.5){$1$}
\htext(2.5 1.5){$0$}
\htext(0.5 0.5){$2$}
\esegment
\move(-3 -15)
\bsegment
\move(1 1)\rlvec(1 0)\rlvec(0 -1)\rlvec(-1 0)\ifill f:0.6
\move(0 0)\rlvec(1 0)\rlvec(0 2)\rlvec(-1 0)\ifill f:0.9
\move(0 2)\rlvec(4 0)
\move(0 1)\rlvec(4 0)\rlvec(0 1)
\move(0 0)\rlvec(2 0)\rlvec(0 2)
\move(3 1)\rlvec(0 1)
\move(1 0)\rlvec(0 2)
\move(0 0)\rlvec(0 2)
\htext(0.5 1.5){$1$}
\htext(1.5 1.5){$1$}
\htext(2.5 1.5){$1$}
\htext(3.5 1.5){$3$}
\htext(0.5 0.5){$2$}
\htext(1.5 0.5){$3$}
\esegment
\move(2.5 -15)
\bsegment
\move(3 2)\rlvec(1 0)\rlvec(0 -1)\rlvec(-1 0)\ifill f:0.6
\move(2 1)\rlvec(1 0)\rlvec(0 1)\rlvec(-1 0)\ifill f:0.9
\move(0 2)\rlvec(5 0)
\move(0 1)\rlvec(5 0)\rlvec(0 1)
\move(0 0)\rlvec(2 0)\rlvec(0 2)
\move(4 1)\rlvec(0 1)
\move(3 1)\rlvec(0 1)
\move(1 0)\rlvec(0 2)
\move(0 0)\rlvec(0 2)
\htext(0.5 1.5){$1$}
\htext(1.5 1.5){$1$}
\htext(2.5 1.5){$1$}
\htext(3.5 1.5){$2$}
\htext(4.5 1.5){$2$}
\htext(0.5 0.5){$2$}
\htext(1.5 0.5){$3$}
\esegment
\move(9 -15)
\bsegment
\move(4 2)\rlvec(1 0)\rlvec(0 -1)\rlvec(-1 0)\ifill f:0.6
\move(3 1)\rlvec(1 0)\rlvec(0 1)\rlvec(-1 0)\ifill f:0.9
\move(1 1)\rlvec(1 0)\rlvec(0 -1)\rlvec(-1 0)\ifill f:0.6
\move(0 0)\rlvec(1 0)\rlvec(0 2)\rlvec(-1 0)\ifill f:0.9
\move(0 2)\rlvec(5 0)
\move(0 1)\rlvec(5 0)\rlvec(0 1)
\move(0 0)\rlvec(3 0)\rlvec(0 2)
\move(4 1)\rlvec(0 1)
\move(2 0)\rlvec(0 2)
\move(1 0)\rlvec(0 2)
\move(0 0)\rlvec(0 2)
\htext(0.5 1.5){$1$}
\htext(1.5 1.5){$1$}
\htext(2.5 1.5){$1$}
\htext(3.5 1.5){$1$}
\htext(4.5 1.5){$2$}
\htext(0.5 0.5){$2$}
\htext(1.5 0.5){$3$}
\htext(2.5 0.5){$3$}
\esegment
\move(15.5 -15)
\bsegment
\move(1 1)\rlvec(1 0)\rlvec(0 -1)\rlvec(-1 0)\ifill f:0.6
\move(0 0)\rlvec(1 0)\rlvec(0 2)\rlvec(-1 0)\ifill f:0.9
\move(0 2)\rlvec(5 0)
\move(0 1)\rlvec(5 0)\rlvec(0 1)
\move(0 0)\rlvec(4 0)\rlvec(0 2)
\move(3 0)\rlvec(0 2)
\move(2 0)\rlvec(0 2)
\move(1 0)\rlvec(0 2)
\move(0 0)\rlvec(0 2)
\htext(0.5 1.5){$1$}
\htext(1.5 1.5){$1$}
\htext(2.5 1.5){$1$}
\htext(3.5 1.5){$1$}
\htext(4.5 1.5){$1$}
\htext(0.5 0.5){$2$}
\htext(1.5 0.5){$3$}
\htext(2.5 0.5){$3$}
\htext(3.5 0.5){$3$}
\esegment
\move(0 0)
\bsegment
\htext(-4.1 -1.3){$1$}
\htext(4.2 -1.3){$2$}
\esegment
\move(0 -2.9)
\bsegment
\htext(-10 -3.5){$1$}
\htext(-3 -3.5){$2$}
\htext(3.3 -3.5){$1$}
\htext(10.2 -3.5){$2$}
\esegment
\move(0 -7.9)
\bsegment
\htext(-14 -3.5){$1$}
\htext(-10.5 -3.5){$2$}
\htext(-5.5 -3.5){$1$}
\htext(-2.4 -3.5){$2$}
\htext(3.6 -3.5){$1$}
\htext(8.2 -3.5){$2$}
\htext(10.7 -3.5){$1$}
\htext(14.3 -3.5){$2$}
\esegment
\move(0 1)
\move(0 -14)
\move(-3.2 -1.1)\ravec(-1 -1)
\move(3.2 -1.1)\ravec(1 -1)
\move(-9 -6.1)\ravec(-1 -1)
\move(-4 -6.1)\ravec(1 -1)
\move(4.2 -6.1)\ravec(-1 -1)
\move(9.2 -6.1)\ravec(1 -1)
\move(-13 -11.1)\ravec(-1 -1)
\move(-10 -11.1)\ravec(0 -1)
\move(-5 -11.1)\ravec(0 -1)
\move(-2 -11.1)\ravec(1 -1)
\move(3 -11.1)\ravec(0 -1)
\move(7 -11.1)\ravec(1 -1)
\move(10.2 -11.1)\ravec(0 -1)
\move(13.2 -11.1)\ravec(1 -1)
\htext(-18 -17){$\vdots$}
\htext(-6 -17){$\vdots$}
\htext(10.5 -17){$\vdots$}
\htext(-11.5 -17){$\vdots$}
\htext(-1.5 -17){$\vdots$}
\htext(4 -17){$\vdots$}
\htext(17 -17){$\vdots$}
\end{texdraw}
\caption{Crystal $\Tinf$ for type $G_2$}\label{fig:3}
\end{figure}

\subsection{Comparison with a previous $A_n$ result}

Our result on $A_n$-type can be found in an earlier work~\cite{hm}.
There, the approach was very different, relying on a work of
Cliff~\cite{cliff}, and the final result was written in a slightly
different form.

The only difference with the current result is that, there, infinitely
many copies of our leftmost column was added to the left of each
representative. This has the advantage of having the Kashiwara
operators look slightly more natural. We do not \emph{insert} or
\emph{remove} columns to remain marginally large, but \emph{push} or
\emph{pull} infinite rows instead.

In the current work, we chose not to add these infinitely many columns,
so as to keep our representatives within the frames of Young tableaux.
The choice between these two presentations seems to be a matter of
taste.

\section{Relationship with another work}

In this section, we recall Cliff's~\cite{cliff} combinatorial
description of $\Binf$ and give an isomorphism between this and our own
realization. Only the finite classical types will be dealt with in this
section, as Cliff did not deal with $G_2$ case.

\subsection{Another realization of $\Binf$}

Let us first recall the crystals
\begin{equation*}
\B_i=\{b_i(k)\vert k\in{\Z}\}
\end{equation*}
defined for each $i\in I$ and introduced in~\cite{MR1240605}. The
crystal $\B_i$ reacts to the Kashiwara operators $\fit$ and $\eit$ by
decrementing or incrementing the inner index $k$, but maps everything
to zero under other $\tilde{f}_j$ and $\tilde{e}_j$. We shall not be
concerned with its exact crystal structure.

Kashiwara has shown~\cite{MR1240605} the existence of an injective
strict crystal morphism
\begin{equation}
\Psi:\Binf\to\Binf\otimes
     \B_{i_k}\otimes\B_{i_{k-1}}\otimes\cdots\otimes\B_{i_1}
\end{equation}
for any sequence $S={i_1,i_2,\cdots,i_k}$ of numbers in $I$ and
Cliff~\cite{cliff} uses this to give a combinatorial description of
$\Binf$ for all finite classical types. To explain this result, we
first define some crystals $\B(i)$ and fix a notation for $\beta_i\in
\B(i)$.
\begin{itemize}
\item
$A_n$ type ($1 \leq i \leq n$)
\begin{equation*}
\B(i)=\B_n\otimes \B_{n-1}\otimes\cdots\otimes \B_i
\end{equation*}
and
\begin{equation*}
\beta_i=b_n(-{k_{i,n}})\otimes b_{n-1}(-{k_{i,n-1}})
\otimes\cdots\otimes b_i(-{k_{i,i}})\in \B(i).
\end{equation*}

\item
$B_n$ type ($1 \leq i \leq n$)
\begin{equation*}
\B(i)=\B_i\otimes\cdots\otimes
      \B_{n-1}\otimes \B_n\otimes \B_{n-1}\otimes\cdots\otimes \B_i
\end{equation*}
and
\begin{equation*}
\beta_i=
\begin{aligned}
&b_i(-{k_{i,\ol{i+1}}}) \otimes b_{i+1}(-{k_{i,\ol{i+2}}})
\otimes\cdots \otimes b_{n-1}(-{k_{i,\bar n}})\\
&\otimes b_n(-{k_{i,n}}) \otimes b_{n-1}(-{k_{i,n-1}}) \otimes \cdots
\otimes b_i(-{k_{i,i}})
\end{aligned}
\in \B(i).
\end{equation*}

\item
$C_n$ type ($1 \leq i \leq n$)
\begin{equation*}
\B(i)=\B_i\otimes \cdots \otimes
      \B_{n-1}\otimes \B_n\otimes \B_{n-1}\otimes\cdots\otimes \B_i
\end{equation*}
and
\begin{equation*}
\beta_i =
\begin{aligned}
&b_i(-{k_{i,\ol{i+1}}}) \otimes b_{i+1}(-k_{i,\ol{i+2}}) \otimes \cdots
\otimes b_{n-1}(-k_{i,\bar n})\\
&\otimes b_n(-k_{i,n})\otimes b_{n-1}(-k_{i,n-1})\otimes \cdots \otimes
b_i(-k_{i,i})
\end{aligned}
\in \B(i).
\end{equation*}

\item
$D_{n+1}$ type ($1\le i\le n-1$)
\begin{align*}
&\B(i)=\B_i \otimes \cdots \otimes \B_{n-1} \otimes \B_{n+1} \otimes
\B_{n} \otimes \cdots \otimes \B_i,\\
&\B(n)=\B_{n+1}\otimes\B_{n},
\end{align*}
and
\begin{align*}
&\beta_i =
\begin{aligned}
&b_i(-k_{i,\ol{i+1}})\otimes
 b_{i+1}(-k_{i,\ol{i+2}})\otimes
 \cdots\otimes b_{n-1}(-k_{i,\ol{n}})\\
 &\otimes b_{n+1}(-k_{i,n+1})\otimes b_{n}(-k_{i,n})\otimes
              \cdots\otimes b_i(-k_{i,i})
\end{aligned}
\in \B(i),\\
&\beta_{n} =
\begin{aligned}
&b_{n+1}(-k_{n,n+1})\otimes  b_{n}(-k_{n,n})
\end{aligned}
\in \B(n).
\end{align*}
\end{itemize}

The following result appears in~\cite{cliff}.
\begin{proposition}
Image of the injective crystal morphism
\begin{equation*}
\Psi : \Binf\rightarrow\Binf\otimes \B(1)\otimes
\B(2)\otimes\cdots\otimes \B(n),
\end{equation*}
is given by
\begin{equation*}
\Psi(\Binf) = \{ u_\infty \otimes \beta_1 \otimes \beta_2 \otimes
\cdots \otimes \beta_n \}.
\end{equation*}
The indices for the components of $\beta_i$ in this expression are
subject to the restrictions given below for each type.
\begin{itemize}
\item $A_n$ type:
\begin{equation*}
0\le k_{i,n}\le k_{i,n-1}\le\cdots\le k_{i,i}
\end{equation*}
for all $1\leq i \leq n$.

\item $B_n$ type:
\begin{equation*}
0\le k_{i,\ol{i+1}}\le k_{i,\ol{i+2}} \le \cdots \le k_{i,\ol n} \le
k_{i,n}/2 \le k_{i,n-1} \le \cdots \le k_{i,i}
\end{equation*}
for all $1\leq i \leq n$. \textup{(}note: $k_{i,n}/2$ need not be an
integer\textup{)}

\item $C_n$ type:
\begin{equation*}
0\le k_{i,\ol{i+1}} \le k_{i,\ol{i+2}} \le \cdots \le k_{i,\ol n} \le
k_{i,n} \le k_{i,n-1} \le \cdots \le k_{i,i}
\end{equation*}
for all $1\leq i \leq n$.

\item $D_{n+1}$ type:
\begin{equation*}
\begin{aligned}
      &0\le k_{i,\ol{i+1}}\le k_{i,\ol{i+2}}\le\cdots\le k_{i,\bar{n}}
      \le\min\bigl(k_{i,n},k_{i,n+1}\bigr)\\
      &\qquad\qquad\qquad\qquad\le \max\bigl(k_{i,n},k_{i,n+1}\bigr)\le
      k_{i,n-1}\le\cdots\le k_{i,i}
\end{aligned}
\end{equation*}
for all $1\leq i \leq n$.
\end{itemize}
\end{proposition}
This proposition gives a combinatorial description of $\Binf$.

\subsection{Isomorphism between two descriptions of $\Binf$}

We now have two explicit descriptions of the same crystal $\Binf$,
namely, $\Psi(\Binf)$ and representatives of $\Tinf$. In this section,
we provide maps between the two in both directions that are crystal
isomorphisms.

With the help of tensor product rules, it is easy to check the
compatibility of these maps with Kashiwara operators. Hence we shall
only write out the maps and give no proofs.

\subsubsection{$A_n$ case}

Element $u_\infty\otimes \beta_1 \otimes \cdots \otimes \beta_n \in
\Psi(\Binf)$ is sent to the tableau whose $i$-th row ($1\le i\le n$)
consists of
      \begin{align*}
      &\text{($k_{i,n}$)-many $(n+1)$\,s,}\\
      &\text{($k_{i,j-1}-k_{i,j}$)-many $j$\,s,
           for each $i<j\le n$, and}\\
      &\text{$\textstyle
      \big((n-i+1)+\sum_{r=i+1}^n k_{r,r}\big)$-many $i$\,s}.
      \end{align*}

\smallskip
Conversely, for each tableau with the $i$-th row consisting of
$b^i_j$-many $j$\,s ($i\prec j\preceq n+1$) and some number of $i$\,s,
we may map it to the element
$u_{\infty}\otimes\beta_1\otimes\cdots\otimes\beta_n \in\Psi(\Binf)$,
where
\begin{equation*}
k_{i,r} = \sum_{j=r+1}^{n+1} b^i_j
\end{equation*}
for $1\le i\le n$ and $i\le r\le n$.

\subsubsection{$B_n$ case}

Element $u_\infty\otimes \beta_1 \otimes \cdots \otimes \beta_n \in
\Psi(\Binf)$ is sent to the tableau whose shape we describe below
row-by-row.
\begin{itemize}\allowdisplaybreaks
\item For $1\le i\le n-1$, the $i$-th row consists of
  \begin{align*}
  &\text{$(k_{i,\ol{i+1}})$-many $\bar i$\,s,}\\
  &\text{$(k_{i,\ol{j+1}}-k_{i,\bar{j}})$-many $\bar j$\,s, for each
        $i<j\le n-1$,}\\
  &\text{$(\lfloor k_{i,n}/2 - k_{i,\bar n}\rfloor)$-many $\bar n$\,s,}\\
  &\text{($(A+B)-(A'+B')$)-many $0$\,s,}\\
  &\text{($\lfloor k_{i,n-1}-k_{i,n}/2\rfloor$)-many $n$\,s,}\\
  &\text{($k_{i,j-1}-k_{i,j}$)-many $j$\,s, for each $i<j\le n-1$, and}\\
  &\text{$\textstyle
  \big((n-i+1)
  + \lfloor \frac{k_{n,n}+1}{2}\rfloor
  + \sum_{r=i+1}^{n-1} k_{r,r}\big)$-many $i$\,s.}
  \end{align*}
\item The $n$-th row consists of
      \begin{align*}
      &\text{$(\lfloor k_{n,n}/2-k_{n,\bar n}\rfloor)$-many
             $\bar n$\,s,}\\
      &\text{($2(B-B')$)-many $0$\,s, and}\\
      &\text{one $n$}.
      \end{align*}
\end{itemize}
Here, $A=k_{i,n-1}-k_{i,n}/2$, $B=k_{i,n}/2-k_{i,\bar n}$, $A'=\lfloor
k_{i,n-1}-k_{i,n}/2]$, and $B'=\lfloor k_{i,n}/2-k_{i,\bar n}\rfloor$,
for each $i\in I$.

\smallskip
Conversely, for each tableau with the $i$-th row consisting of
$b^i_j$-many $j$\,s ($i\prec j\preceq \ol{i}$) and some number of
$i$\,s, we may map it to the element
$u_{\infty}\otimes\beta_1\otimes\cdots\otimes\beta_n \in\Psi(\Binf)$,
where, for $1\le i\le n-1$,
\begin{alignat*}{2}
&k_{i,r}
      =\sum_{j=r+1}^{n} b^i_j+b_0^i+\sum_{j=i}^{n} b^i_{\bar{j}}
      & &\qquad\text{for $i\le r\le n-1$},\\
&k_{i,n}
      =2(\sum_{j=i}^{n} b^i_{\bar{j}})+b_0^i,\\
&k_{i,\bar{r}}
      =\sum_{j=i}^{r-1} b^i_{\bar{j}}
      & &\qquad\text{for $i+1\le r\le n$},
\end{alignat*}
and
\begin{equation*}
\text{$k_{n,n}=2b^n_{\bar{n}}+b_0^n$}.
\end{equation*}

\subsubsection{$C_n$ case}

Element $u_\infty\otimes \beta_1 \otimes \cdots \otimes \beta_n \in
\Psi(\Binf)$ is sent to the tableau whose shape we describe below
row-by-row.
\begin{itemize}\allowdisplaybreaks
\item For $1\le i\le n-1$, the $i$-th row consists of
      \begin{align*}
      &\text{$k_{i,\ol{i+1}}$-many $\bar i$\,s,}\\
      &\text{$(k_{i,\ol{j+1}}-k_{i,\bar{j}})$-many $\bar j$\,s,
           for each $i<j\le n-1$},\\
      &\text{$(k_{i,n}-k_{i,\bar{n}})$-many $\bar n$\,s},\\
      &\text{($k_{i,j-1}-k_{i,j}$)-many $j$\,s,
           for each $i<j\le n$, and}\\
      &\text{$\textstyle
      \big((n-i+1)+\sum_{r=i+1}^n k_{r,r}\big)$-many $i$\,s}.
      \end{align*}
\item The $n$-th row consists of
      \begin{align*}
      &\text{$k_{n,n}$-many $\bar n$\,s and}\\
      &\text{one $n$}.
      \end{align*}
\end{itemize}

\smallskip
Conversely, for each tableau with the $i$-th row consisting of
$b^i_j$-many $j$\,s ($i\prec j\preceq \ol{i}$) and some number of
$i$\,s, we may map it to the element
$u_{\infty}\otimes\beta_1\otimes\cdots\otimes\beta_n \in\Psi(\Binf)$,
where, for $1\le i\le n-1$,
\begin{alignat*}{2}
&k_{i,r}
      =\sum_{j=r+1}^{n} b^i_j+\sum_{j=i}^{n} b^i_{\bar{j}}
      & &\qquad\text{for $i\le r\le n-1$},\\
&k_{i,n}
      =\sum_{j=i}^{n} b^i_{\bar{j}},\\
&k_{i,\bar{r}}
      =\sum_{j=i}^{r-1} b^i_{\bar{j}}
      & &\qquad\text{for $i+1\le r\le n$},
\end{alignat*}
and
\begin{equation*}
\text{$k_{n,n}=b^n_{\bar{n}}$},
\end{equation*}

\subsubsection{$D_{n+1}$ case}

Element $u_\infty\otimes \beta_1 \otimes \cdots \otimes \beta_n \in
\Psi(\Binf)$ is sent to the tableau whose shape we describe below
row-by-row.
\begin{itemize}\allowdisplaybreaks
\item For $1\le i\le n-1$, the $i$-th row consists of
      \begin{align*}
      &\text{$k_{i,\ol{i+1}}$-many $\bar i$\,s,}\\
      &\text{$(k_{i,\ol{j+1}}-k_{i,\bar{j}})$-many $\bar j$\,s,
      for each $i<j\le n-1$},\\
      &\text{$(k_{i,n}-k_{i,\bar{n}}-
           \text{max}\{0,k_{i,n}-k_{i,n+1}\}$}\\
      &\text{$=k_{i,n+1}-k_{i,\bar{n}}-
           \text{max}\{0,k_{i,n+1}-k_{i,n}\})$-many
           $\bar{n}$\,s},\\
      &\text{$(\text{max}\{0,k_{i,n+1}-k_{i,n}\})$-many
           $(\ol{n+1})$\,s},\\
      &\text{($\text{max}\{0,k_{i,n}-k_{i,n+1}\}$)-many
           $(n+1)$\,s},\\
      &\text{$(k_{i,n-1}-k_{i,n}-
           \text{max}\{0,k_{i,n+1}-k_{i,n}\}$}\\
      &\text{$=k_{i,n-1}-k_{i,n+1}-
           \text{max}\{0,k_{i,n}-k_{i,n+1}\})$-many
           $n$\,s},\\
      &\text{($k_{i,j-1}-k_{i,j}$)-many $j$\,s,
           for each $i<j\le n-1$, and}\\
      &\text{$\textstyle \big(
      C+(n-i+1)+ \sum_{r=i+1}^{n-1}k_{r,r}\big)$-many $i$\,s}.
      \end{align*}
\item The $n$-th row consists of
      \begin{align*}
      &\text{$(k_{n,n}-
           \text{max}\{0,k_{n,n}-k_{n,n+1}\}$}\\
      &\text{$\quad =k_{n,n+1}-
           \text{max}\{0,k_{n,n+1}-k_{n,n}\})$-many
           $\ol{n}$\,s,}\\
      &\text{$(\text{max}\{0,k_{n,n+1}-k_{n,n}\})$-many
           $(\ol{n\!+\!1})$\,s},\\
      &\text{($\text{max}\{0,k_{n,n}-k_{n,n+1}\}$)-many
           $(n\!+\!1)$\,s, and}\\
      &\text{one $n$}.
      \end{align*}
\end{itemize}
Here, $C = k_{n,n}+\max\{0,k_{n,n+1}-k_{n,n}\} = k_{n,n+1} + \max\{0,
k_{n,n}-k_{n,n+1}\}$.

Conversely, for each tableau with the $i$-th row consisting of
$b^i_j$-many $j$\,s ($i\prec j\preceq \ol{i}$) and some number of
$i$\,s, we may map it to the element
$u_{\infty}\otimes\beta_1\otimes\cdots\otimes\beta_n \in\Psi(\Binf)$,
where, for $1\le i\le n-1$,
\begin{alignat*}{2}
&k_{i,r}
      =\sum_{j=r+1}^{n+1} b^i_j+\sum_{j=i}^{n+1} b^i_{\bar{j}}
      &\qquad&\text{for $i\le r\le n-1$},\\
&k_{i,n}
      =b_{n+1}^i+\sum_{j=i}^{n} b^i_{\bar{j}},\\
&k_{i,n+1}
      =\sum_{j=i}^{n+1} b^i_{\bar{j}},\\
&k_{i,\bar{r}}
      =\sum_{j=i}^{r-1} b^i_{\bar{j}}
      &&\text{for $i+1\le r\le n$},
\end{alignat*}
and
\begin{align*}
&\text{$k_{n,n}
      =b_{n+1}^n+\sum_{j=n}^{n} b^n_{\bar{j}}$},\\
&\text{$k_{n,n+1}
      =\sum_{j=n}^{n+1} b^n_{\bar{j}}$}.
\end{align*}


\end{document}
